\documentclass[11pt]{article}
\usepackage{float}
\usepackage[margin=1in]{geometry}
\usepackage{amsmath,amssymb,amsthm,mathtools}
\usepackage{enumitem}
\usepackage[hidelinks,hypertexnames=false]{hyperref}
\usepackage[numbers]{natbib}
\newtheorem{theorem}{Theorem}[section]
\newtheorem{proposition}[theorem]{Proposition}
\newtheorem{lemma}[theorem]{Lemma}
\theoremstyle{definition}
\newtheorem{assumption}[theorem]{Assumption}
\newtheorem{definition}[theorem]{Definition}

\theoremstyle{remark}
\newtheorem{remark}[theorem]{Remark}
\newtheorem*{remark*}{Remark}
\newcommand{\vs}{\vskip.05in}
\newcommand{\R}{\mathbb R}
\newcommand{\N}{\mathbb N}
\newcommand{\E}{\mathbb E}
\newcommand{\Pp}{\mathbb P}
\newcommand{\1}{\mathbf 1}

\newcommand{\cM}{\mathcal M}
\newcommand{\cE}{\mathcal E}
\newcommand{\cS}{\mathcal S}
\newcommand{\eps}{\varepsilon}
\newcommand{\ep}{\varepsilon}
\newcommand{\weakto}{\rightharpoonup}
\newcommand{\dd}{\,\mathrm d}
\newcommand{\be}{\begin{equation}}
\newcommand{\eeq}{\end{equation}}
\newcommand{\ee}{\end{equation}}

\newcommand{\ol}{\overline}
\newcommand{\beq}{\begin{equation}}
\newcommand{\pt}{\partial}
\newcommand{\wh}{\widehat}
\newcommand{\oo}{\infty}

\usepackage{parskip}
\usepackage[dvipsnames,table,xcdraw]{xcolor}
\title{The microstructure of information through mean field learning
interaction}

\author{ Panagiotis E. Souganidis \thanks {Department of Mathematics, The University of Chicago.  Partially supported by the NSF grants DMS-2153822 and DMS-2452972. \textbf{Email:} souganidis@math.uchicago.edu} \and  Renyuan Xu \thanks{Department of Management Science and Engineering, Stanford University. Partially supported  by the NSF CAREER Award DMS-2524465 and NSF grant DMS-2602037. \textbf{Email:} renyuanxu@stanford.edu} \and Thaleia Zariphopoulou \thanks{Departments of Mathematics and IROM, The University of Texas at Austin, Austin,
USA,  and the Oxford-Man Institute, University of Oxford, Oxford, UK. \textbf{Email:} zariphop@math.utexas.edu}}
\date{\today}

\begin{document}
\maketitle
\begin{abstract}
\noindent We study how knowledge about partially observable model parameters is being created and propagated in large scale systems. Information is generated from observational signals whose dynamics depend endogenously on the learning history, the learning duration, the decisions and the related costs encountered by a continuum of decision makers. Such characteristics give rise to a new class of sequential testing problems in which there is continuous endogenous feedback. Adopting the classical approach of {Shiryaev}, we relate the problem to an optimal stopping one for the representative belief process, whose generator now depends on the evolving interactions. Considering a fixed point for the latter, we introduce a notion of mean field game and prove the existence of an equilibrium. Our arguments, which  are based on a combination of pde and probabilistic techniques, involve a detailed analysis of the underlying variational inequality and the geometry, stability and continuity of the free boundaries. % using a combination of pde and probabilistic techniques.    
\vs
\vs

\end{abstract}

\section{Introduction}

This paper studies how interactive and collective information acquisition about unknown model parameters generates and shapes aggregate learning about them in large-scale systems. 
Information is dynamically created in a Bayesian framework by costly signals which generate belief processes (posterior, conditional on the signals, distributions). 
The key feature in our work is that the stochastic dynamics of the signals are directly affected by the collective information-seeking behavior, like duration of learning and  
assessment about the unknown parameter and decisions when learning terminates.  Acquiring information induces a cost consisting of a running component while signals are being used and a final aggregate cost caused by these possibly erroneous assessments and decisions at exit. The focus is on how to (i)~minimize the representative information acquisition cost, (ii)~ investigate through its solution how representative knowledge is being created and shaped by a continuous interaction of signals, (iii)~examine its evolution and (iv)~build an equilibrium model for its validity. 
Our results spread across the areas of hypothesis testing, Bayesian learning, filtering, optimal stopping, and mean field games with 
partial observation.
%Sequential decision-making problems under uncertainty arise in many
%large-scale systems where individual agents lack direct access to the true
%underlying state. In such settings, each agent must construct beliefs from
%noisy observations by acquiring information over time and subsequently
%determine optimal control actions, such as stopping times, based on these
%beliefs. Classical models view information flows as exogenous and examine
%optimal strategies under state-dependent stochastic signals.
\vs
The motivation to study large-scale models of interactive learning comes from numerous applications in marketing, advertising, trading, autonomous  driving systems, and supply-chain systems, to name a few. Indeed, in e-commerce platforms, knowledge about the quality of a new product is generated from reviews and usage data induced by earlier buyers; see Kakhbod,
Lanzani and Xing \cite{kakhbod2021heterogeneous}. Because some consumers
provide detailed feedback while others leave only short comments, the
informativeness of the public, collective signals {depends} considerably on how many and who has chosen
to adopt the new product so far. Indeed, if early buyers provide meaningful, detailed feedback, then
the signal about the quality of the product is strong. On the other hand, if many early buyers provide only
sparse or uninformative feedback, the collective signal is inferior. As a result,
the quality of acquired information is being shaped directly by
the previous adoption decisions, the population size of this group and the duration of their learning.
\vs

In financial markets with liquidity effects, if asset liquidation becomes more intense, the overall market depth
deteriorates, making price movements more volatile and less informative
about the underlying price fundamentals. As a consequence, the ability to
extract information from observed prices depends critically on how many positions are already closed and by when this actually occurred; see, for example, Casgrain and Jaimungal 
\cite{casgrain2020mean}. The resulting knowledge acquisition problem is thus affected
critically by collective learning interaction, as the value of waiting to gather
further information and the quality of information itself are determined by
the evolving distribution of decisions and learning duration.

\vs

The case of a single decision maker, who learns about a binary random variable $\Theta $ with possible values $0$ and $1$, has been extensively analyzed in the literature, starting
with the  work of Wald in \cite{wald1947foundations}\ on costly sequential learning in Bayesian statistics. Since then, it has a long history of
development and constitutes a central area in statistical learning and hypothesis testing, with hypotheses being $H_{0}:\Theta =0$ and $H_{1}:\Theta =1$, for this binary case.  The underlying problem is to minimize, up to the time that learning stops, the total expected cost {consisting} of a running part due to signal usage and a terminal component. The latter is generated by erroneous assessments about the values of $\Theta $ at exit time and related decisions. Thus, the controls are formulated in terms of the pair of stopping times and decisions at them, $(\tau(\omega); d_i(\tau(\omega)))$ for $i=0,1$.
\vs
{Shiryaev} developed in  \cite{shiryaev2007optimal} a connection to an optimal stopping problem for the belief process. In particular, it was shown in \cite{shiryaev2007optimal} that the underlying free boundary problem has two continuous and not intersecting boundaries which completely characterize the optimal learning. Specifically, the boundaries divide the space in an information seeking region and two idle regions to its left and right sides. If the belief process starts, at $t=0$ in the left (resp. right) region, no learning occurs and hypothesis $H_{0}:\Theta =0$ (resp. {$H_{1}:\Theta =1$}) is instantly accepted. The cost then comes from the instant decision $d_0(0)$ (resp. $d_1(0)$). If, on the other hand, the belief process starts in the middle region, learning takes place till the first hitting time to either {boundary}. If the left boundary (resp. right) is reached first  at $\tau^*(\omega)$, hypothesis $H_{0}:\Theta =0$ (resp. {$H_{1}:\Theta =1$}) is accepted and decision $d_0(\tau^*(\omega))$ (resp. $d_1(\tau^*(\omega)$) is made. The total cost is the one {incurred} by the latter and the running learning cost in $[0,\tau^*(\omega)]$. 
%We will base our solution approach on this  connection to analyze the mean field model we propose herein.
\vs
The results of  \cite{shiryaev2007optimal} initiated a very rich and active line of research, and providing a complete bibliography is beyond the scope herein (see, among others, Moscarini and Smith \cite{moscarini2001optimal}, Miao and Xing \cite{miao2024dynamic}, Ke, Tang, Villas-Boas and Zhang \cite{ke2022parallel}). We also mention the works of Ekstr\"om, Karatzas and Vaicenavicius \cite{EKV}, Ekstr\"om and Vaicenavicius \cite{EkVa}, Ekstr\"om and Wang \cite{EkWa}  and Campbell, Gaitsgori, Groenewald and Karatzas \cite{CaGaGrKa}. 
\vs

Herein, we build a mean field game of interactive and collective information acquisition, which we study using the technical approach of \cite{shiryaev2007optimal}. We assume that information is being acquired and propagated through a continuous interaction within a population continuum which uses signals with dynamics directly affected by the aggregate learning, decisions and costs. This continuum starts with a distribution of prior probabilities for the values of $\Theta$ and interacts through time till learning stops. The key quantity is the temporal aggregate of the fraction of the population that has completed learning, selected a decision, encountered the relevant costs and exited. This temporal aggregate alters the learning mechanism of the remaining population through the drift of its, continuously updated, signal. In turn, it affects the representative value function which models the minimal total expected informational cost and constitutes the best response to the effects of realized learning and costs. We are then motivated to introduce a fixed-point notion in terms of this aggregate feedback which leads to a mean field game. 
% we propose herein.
\vs
We stress that contrary to the standard mean field games where players are entirely identified by their controlled state processes, herein the notion of a player is identified with the individual prior probability and subsequent belief process. In other words, there is no representative controlled diffusion as such but rather a representative belief process whose dynamics are affected by the learning and decisions of the peer continuum. To our knowledge, the only work close to ours is the preprint of Campbell and Zhang  \cite{campbell2024mean}, but both their technical assumptions and modeling elements are different than ours. We provide further details at the end of the introduction. 
\paragraph{The model and the results.} In what follows we explain our main results and  refer for the precise definitions and assumptions later in the paper. 
\vs
 Let $(\Omega,\mathcal F,\Pp)$ be a probability space that supports 
 \beq\label{hob1}
\begin{split}
&\text{(i)~a binary random variable  $\Theta \in \{0,1\}$,}\\
&\text{(ii)~a random variable $P_0\in [0,1]$ with {law $\eta$ having a continuous density}}\\
&\quad \text{ representing the initial belief (prior) that $\Pp(\Theta=1{\mid P_0})=P_0$,}\\
&\text{(iii)~a standard Brownian motion $W$.}
\end{split}
\eeq
It is assumed that 
\beq\label{hob2}
\text{{$W$ is independent of $(\Theta,P_0)$.}}
\eeq
 
% an initial belief $P_0$ with law with the property that
%\beq\label{density}
%\eta \ \ \text{is   probability measure 
%on $[0,1]$ with density  denoted also by $\eta$,} 
%\eeq
%%taken to be  probability measure 
%%on $[0,1]$ with density which we denote also by $\eta$, 
%an
%unobserved binary state $\Theta\in\{0,1\}$, and a Brownian motion $W$.
%We assume % that $P_0$ has law $\eta$, a probability measure on $[0,1]$,
%that 
%\beq\label{vlakas20}
%\Pp(\Theta=1\mid P_0)=P_0 \ \ \text{and} \ \ \text{$W$ is independent of $(P_0,\Theta)$},
%\eeq
%and we 
%denote by $\mathcal M_+([0,T])$ and $\mathcal P([0,T])$ the {nonnegative Borel measures with total mass at most one}
%and  the probability measures on $[0,T]$ respectively.  
%\vs

We study information acquisition about the  unknown binary state $\Theta$ from continuous time signals in a large-scale system. Knowledge is being created from a representative observational signal $Y^\mu$ whose dynamics depend on aggregate learning interaction in the continuum. 
Specifically, $Y^\mu$ is taken to solve 
\begin{equation}\label{vlakas0}
\dd Y^\mu_t =b\big(t,\Theta;F(t;\mu_{\le t})\big) \dd t + \dd W_t  \ \ \text{for} \ \    t\in [0,T]  \ \ \text{and} \ \ Y^\mu_0=0.
\end{equation}
The term $F(t;\mu_{\le t})$, which is defined  in \eqref{eq:causal-smoothing}, models the collective information feedback and incorporates the learning duration, exit times, and decisions about the unknown values and related costs  up to time $t$. The input  $\mu$  is a pair $(\mu^0,\mu^1)$ of two nonnegative subadditive measures  representing stopping times and decisions at them,  given, for each Borel subset $A$ of $[0,T]$ and $i=0,1$, by
\begin{equation*}
\mu^i(A)=\mathbb{P}(\tau \in A,d=i), % \, \, \text{ for }\, \, i\in \{0,1\}.
%A\in \mathcal{B}%
%([0,T]) \, \, \text{ and }\, \, i\in \{0,1\},
\end{equation*}
and such that $\mu^0+\mu^1$ is a probability measure on $[0,T]$. In what follows, we denote the set of such pairs by $\cM$, that is,
\begin{equation}
 \mathcal M:=\Bigl\{(\mu^0,\mu^1):
 \mu^0,\mu^1\in\mathcal M_+([0,T]), {\mu^0+\mu^1\in \mathcal{P}([0,T])}\Bigr\},
 \label{eq:set_C}
\nonumber
\end{equation}
where $\mathcal M_+([0,T])$ stands for the set of nonnegative subadditive measures on $[0,T]$.
\vs
For a fixed input $\mu$, the observation signal  $Y^\mu$ generates the representative posterior belief process $P$ given by 
%(P_t)_{0\leq t\le T}$ 
%given by 
\begin{equation*}
P_t:=\mathbb{P}(\Theta=1\mid {\mathcal{F}_t^\mu}) \ \ \text{for} \ \  t\in [0,T] \ \ \text{and} \ \ P_0=p\in [0,1],
\end{equation*}
where  $\mathbb F^\mu=(\mathcal F_t^\mu)_{0\le t\le T}$ denotes the usual
augmentation of the filtration
$(\sigma(P_0,Y_s^\mu) :0\le s\le t))_{0\le t\le T}$.
%and the initial belief $P_0$ is observed by the agent.
\vs
Classical  results in filtering yield that the process $P$ solves the sde
\begin{equation}\label{0}
\, \mathrm{d }P^\mu_t =\Delta b^\mu(t)P^\mu_t(1-P^\mu_t)\, \mathrm{d }Z^\mu_t \ \ \text {for } \ \ t\in (0,T] \ \ \text{and} \ \ P_0=p\in [0,1], 
\end{equation}
where 
\be\label{ze10}
\Delta b^\mu(t):=\Delta b\bigl(t;F(t;\mu_{\le t})\bigr)={b(t,1; F(t;\mu_{\le t}))-b(t,0; F(t;\mu_{\le t}))}
\ee
is the difference between the signal drifts under the two
states $0$ and $1$, and $Z^\mu$ is a $\mathbb F^\mu$-Brownian motion referred to as the innovation process; see, for example, Bain and Crisan \cite{BainCrisan2009}. 
\vs
We assume (see \eqref{takis11} and \eqref{ass.b}) that $b$ and $\Delta b$ satisfy the assumptions that are needed for 
\eqref{vlakas0} and \eqref{0} to have unique strong solutions. 
%\vs
%{\bf 
%We denote by $\overline G(p,0)$  and $\overline G(p,1)$ the cost incurred by choosing state $0$ and, respectively, $1$ at posterior belief $p$.
%%$\overline G(p,1)$ for state $1$, at posterior belief $p$.
%The loss at stopping after minimizing over the
%decision is assumed  to be 
%\beq\label{0.0}
%G(p)=\min \{\overline G(p,0), \overline G(p,1)\}.
%\eeq}
%\vs
%We also assume that during learning there is a  running cost
%\beq\label{00}
%f=f(t,p; \mu).
%\eeq
%%{\bf During learning, a running cost also occurs which -  Takis: we need to say how from a general costs in the original formulation, we end up with a Makrovian running cost $f(t,P)$. Perhaps, I should go back and add something in the paragraph for Shirayev? or do we just add a sentence related to this? }
%
\vs
The   representative information acquisition problem is an optimization one {over} all stopping learning times 
and decisions related to parameter value assessment.
\vs
Following  \cite{shiryaev2007optimal}, we formulate the original problem as an optimal stopping problem for the representative belief process over all stopping learning times.
% and decisions relasted to parameter value assessment.
\vs
Specifically,  we introduce the value function 
\begin{equation}\label{000}
V^\mu(0,p) =\inf_{\tau  \in \mathcal{T}^\mu_{0,T}} \mathbb{E}\left[ \int_0^{\tau}
f^\mu(s,P^\mu_s)\, \mathrm{d }s+G(P^\mu_{\tau}) :\, P^\mu_0=p \right],
\end{equation}
where, for $t\in [0,T]$, $\mathcal{T}^\mu_{t,T}$ denotes the set of $[t,T]$-valued
{$\mathbb F^\mu$-stopping times}, $f^\mu: [0,T]\times [0,1] \to \R$ is the running cost during learning given by 
\be\label{eq:f-intro}
f^\mu(t,p):=f(t,p;F(t;\mu_{\le t})),
 %:=\ol f\bigl(t,p;F(t;\mu_{\le t})\bigr),
 \eeq
with $f$ satisfying \eqref{eq:primitive-bounds-1}, and 
$G:[0,1]\to \R$ is {the penalty incurred when learning stops}, taken to be
\beq\label{chi1}
G(p):=\underset{i=0,1}\min \, \ol G(p,i).
\nonumber
\eeq
Following \cite{shiryaev2007optimal} and subsequent papers, we assume that, for  $a_0,a_1>0$, 
\beq\label{chi2}
\ol G(p,0)= a_0 p \ \ \text{and} \  {\ol G}(p,1)=a_1(1-p).
\nonumber
\eeq
With this choice of $\ol G$, the terminal cost $G$ can be written as  
\beq\label{chi4}
G(p)=\min\{a_0 p, a_1(1-p)\}.
\eeq
The maximum in \eqref{chi4} is attained at 
\be\label{ze2.1}
m=\frac{a_1}{a_0+a_1}.
\ee
Then, we can rewrite $G$ as
\beq\label{chi3}
G(p)={\overline G}(p,D(p)),
\nonumber
\eeq
where $D:[0,1]\to \{0,1\}$ is the so-called decision rule, with the tie at $m$ assigned to decision $0$, 
\begin{equation}
 D(p):=\begin{cases}
 0& \textrm{ if }  \,\, p\le m,\\
 1& \textrm{ if } \,\, p>m.
 \end{cases}
 \label{eq:optimal_d}
\nonumber
\end{equation}
We remark that it is possible to work with  other choices of $\ol G(\cdot,0)$ and $\ol G(\cdot,1)$ as long as each one  of them is concave and $\ol G(0,0)<\ol G(0,1)$ and {$\ol G(1,0)>\ol G(1,1)$}, provided we make assumptions on the sign 
of $\pt_p^2 G$ near the maximum. In addition, we can assume that the constants $a_0,a_1$ may depend on time. Here,  for the {sake} of simplicity and in order not to confuse the reader with unnecessary mathematical technicalities, we have chosen to work with the classical choice of $G$ in \eqref{chi4}.
\vs 

The population learning feedback $\mu \in \cM$ influences not only the dynamics of $P^\mu$ but, also,  the running cost, and $V^\mu$ is interpreted as a best response towards $\mu$ of the knowledge creation problem. This leads to a notion of mean field equilibrium, which {we} formulate in the sequel.
%by looking at a fixed point of a map $\mathcal{M} \ni \mu\longrightarrow\Phi(\mu) \in \mathcal{M}$
%defined by \eqref{eq:fixed point}.
\vs

For each fixed $\mu \in \cM$, the time dependent value function $V^\mu : [0,T]\times [0,1]\to \R$ is defined, in complete analogy to \eqref{000}, by
\beq\label{hob5}
V^\mu(t,p) =\inf_{\tau  \in \mathcal{T}^\mu_{t,T}} \mathbb{E}\left[ \int_t^{\tau}
f^\mu(s,P^{t,p;\mu}_s)\, \mathrm{d }s+G(P^{t,p;\mu}_{\tau}) \right],
\nonumber
\eeq
where $P^{t,p;\mu}$ solves \eqref{0} in $[t,T]$ and $P^{t,p;\mu}_t=p$.
\vs
We establish that  $V^\mu$ 
is the unique (viscosity) solution to  the terminal value variational inequality 
\begin{equation}\label{takis100}
\max \Big \{-\pt_t v -a^\mu \pt^2_p v-f^\mu,\,v-G \Big \}=0 \ \ \text{in} \ \ [0,T)\times [0,1] \qquad v(T,p)={G(p)},
\end{equation}
where  
\begin{equation}\label{ze1}
a^\mu(t,p) =\frac12\bigl(\Delta b^\mu(t\bigr)\big)^2p^2(1-p)^2.
%\frac12\bigl(\Delta b(t;\mu_{\le t})\bigr)^2p^2(1-p)^2.
\end{equation}
%%and the terminal value $v(T,p)=G(p)$. 
%\vs 
%We assume that $\Delta b^\mu$ is uniformly positive, bounded and Lipschitz, and $f^\mu$ and $G$ are   
%bounded, Lipschitz and concave. 
%The precise assumptions on all the coefficients  in  \eqref{takis100} are stated in section~2. %
%$f^\mu(t,\cdot)$ and $G$ concave; the precise assumptions are stated in section~2. 
% while $G$ is concave and either piecewise affine or
%smooth; see Assumptions~\ref{ass:coefficients}--\ref{ass:obstacle}. 

%We establish that $V^\mu(t,p)$ is the unique continuous state-constraint viscosity solution of (). 

We show that, for each $\mu \in \cM$,  the domain $[0,T]\times [0,1]$ is divided in two regions, namely,  the 
exploration region 
\be\label{vlakas01}
\cE^\mu=\{(t,p) \in [0,T]\times [0,1]:\, V(t,p)< G(p)\},
\eeq
and the stopping region 
\beq\label{vlakas02}
\cS^\mu=\{(t,p) \in [0,T]\times [0,1]:  V(t,p)= G(p)\}.\eeq
In what follows, for each $t \in [0,T]$, we write $\cE^\mu_t$ and $\cS^\mu_t$ for the $t-$slices of the exploration and stopping regions respectively, and note that $\cE^\mu_T=\emptyset$ and $\cS^\mu_T=[0,1].$ We  show that 
\begin{equation*}
\mathcal{E}_t^\mu=(A_-^\mu(t),A_+^\mu(t)) \ \ \text{and} \ \ \cS^\mu_t=[0,A_-^\mu(t)] \cup [A_+^\mu(t),1], 
\end{equation*}
%We also establish that, for each $t\in [0,T)$,
% \beq\label{zoe1}
%\cE^\mu_t=(A_-^\mu(t), A_+^\mu(t)) \ \ \text{and} \ \ \cS^\mu_t=[0,A_-^\mu(t)] \cup [ A_+^\mu(t),	1], 
%\eeq
where, for $t\in [0,T]$, 
%$\cE_t^\mu$ and $\cS_t^\mu$ are the $t-$ slices of $\cE^\mu$ and $\cS^\mu$ and 
\begin{equation}
 A_-^\mu(t):=\sup\big\{p \in \cS_t^\mu  \cap [0,m]\big\} \ \ \text{and} \ \ 
 A_+^\mu(t):=\inf \big \{ p \in \cS_t^\mu \cap(m,1]\big\}.
 \label{eq:boundary-definitions-intro}
\end{equation}
The maps $t\to A_-^\mu(t)$ and $t\to A_+^\mu(t)$ are, respectively,  the left and right free boundaries of $V^\mu$. 
\vs
We show that the value function is uniformly stable with respect to the weak convergence of $\cM$, that is,  if $\mu^n, \mu \in \cM$ are such that $\mu^n \weakto \mu$, then 
the  value functions $V^{\mu^n}$ converge uniformly on $[0,T]\times [0,1]$ to $V^\mu$. 
\vs
%are also uniformly continuous with respect to the weak convergence of measures in $\cM$. 
In turn, we combine this continuity with various local regularity and
non-degeneracy estimates to prove that there exists $\mathcal R\subset (0,T)$,  consisting of the times where the interfaces  $[0,T]\ni t\to A_\pm^\mu(t)$ are continuous,  of full measure such that, as $\mu^n \rightharpoonup \mu$ and locally uniformly in  $\mathcal R$,
%\mu^n\rightharpoonup \mu$, Lemma~\ref{lem:coefficient-stability} and
%Proposition~\ref{prop:value-stability} yield 
%\begin{equation*}
%\|V^{\mu^n}-V^\mu \|_\infty \longrightarrow0.
%\end{equation*}
%Combining the above stability result with various local regularity and
%non-degeneracy estimates, together with the regularity condition in \cite[%
%Theorem~1]{BDM2010}, we prove 
\begin{equation}\label{vlakas04}
A_\pm^{\mu^n}\longrightarrow A_\pm^\mu.
\nonumber
\eeq
%The set  $\mathcal R$ consists of the times where the interfaces  $A_\pm^\mu$ are continuous. 
%\vs
%\ \ \text{locally uniformly on} \ \ (0,T) \ \ \text{and}  \ \  A_\pm^{\mu^n}(0)\longrightarrow 
%A_\pm^\mu(0).
%\end{equation}
In the sequel, we use  the best response pair $(\hat\tau_{0,p}^\mu, \hat d_{0,p}^\mu)$, which, for  each  initial belief $p$ and  information input $\mu\in \cM$, %we use the best response pair $(\hat\tau_{0,p}^\mu, \hat d_{0,p}^\mu)$
is given by
\begin{eqnarray}\label{eq:best-response-pair}
    \hat\tau_{0,p}^\mu
 :=\inf\{t\in[0,T]:
          V^\mu(t,P_t^{0,p;\mu})=G(P_t^{0,p;\mu})\} \ \ 
\textrm{ and } \ \ 
 \hat d_{0,p}^\mu
 :=D(P_{\hat\tau_{0,p}^\mu}^{0,p;\mu}). 
\end{eqnarray}
%where $\inf\varnothing:=T$.  
%which is, as we establish later, %in view of the dynamic programming principle, 
%is optimal. 
% (see section \ref{sec:VI-boundaries}).  Since the posterior solution map \((p,\omega,t)\mapsto P_t^{0,p;\mu}(\omega)\) is jointly measurable and continuous in \(t\), and \(V^\mu\) is continuous, the best response pair is measurable in \((p,\omega)%\).  
%\vs

For each  Borel set
$A\subset[0,T]$ and  $i=0,1,$ we define 
\[
 \nu_i^\mu(A,p)
 :=\Pp(\hat\tau_{0,p}^\mu\in A,\ \hat d_{0,p}^\mu=i) \ \  \text{and} \ \ 
 \nu_i^\mu(A):=\int_{[0,1]}\nu_i^\mu(A,p)\,\dd \eta(p),
\]
%\begin{remark}%[Initial belief distribution $\eta$]
%The initial belief being a random variable with distribution $\eta$ models the heterogeneity of prior knowledge in the large scale system. This differs from the setting of \cite{campbell2024mean}, where it is assumed that there only common  prior belief $p_0$.
%\end{remark}
and note that $(\nu_0,\nu_1) \in \c M(0,T)$. 
%It is that $ that $ c_0^\mu([0,T])+ c_1^\mu([0,T])=1.$
\vs 
Next, we introduce  the map 
\begin{equation}
 \cM\ni \mu \to \Phi(\mu):=(\nu_0^\mu,\nu_1^\mu)\in\mathcal M,
 \label{eq:fixed point}
\nonumber
\end{equation}
and define a mean field (solution) equilibrium  to be  a fixed point
\begin{equation}
 \Phi(\mu^*)=\mu^*.
 \label{eq:fixed-point-general}
\nonumber
\end{equation}
Equivalently, for $i=0,1$ and any Borel $A\subset[0,T]$, the components $(\mu^{0,\ast}(A), \mu^{1,\ast}(A))$ of $
\mu^\ast$ satisfy the consistency relationship
\[
 \mu^{*,i}(A)
 =\int_{[0,1]}
   \Pp(\hat\tau_{0,p}^{\mu^*}\in A,
       \hat d_{0,p}^{\mu^*}=i)\,\dd \eta(p).
\]
%To prove existence, it remains to establish continuity of $\Phi$. Section~3 develops the value-function and stopping-boundary stability needed for this step; Section~4 then applies the Schauder--Tychonoff theorem.

The main result of the paper is stated next.

\begin{theorem}\label{main}
Assume \eqref{takis20}. % and either \eqref{g1} or \eqref{g2}. 
Then, the map $\mathcal M \ni \mu \to \Phi(\mu)$ has a fixed point. 
\end{theorem}

%locally uniformly in $(0,T)$. 
For each candidate population learning distribution $\mu \in \mathcal{M}$, the joint law of the best response defines $\Phi(\mu)\in \cM$. We show that the map $\Phi:\cM\to \cM$ admits a fixed point $\mu^\ast$, which represents {a} 
mean field game equilibrium. 

%joint law of the best response defines $\Phi(\mu)\in \mathcal{C}$. We show
%that $\Phi$ admits a fixed point $\Phi(\mu^*)=\mu^*$. This fixed-point
%represent a mean field game equilibrium. 
%A fixed point $\Phi(\mu^*)=\mu^*$ defines exactly a mean field equilibrium. Our stability result ensures the continuity of $\Phi$  (Proposition~\ref{prop:canonical-response-continuity}). Since $\mathcal C$ is compact and convex, the Schauder--Tychonoff theorem yields such a fixed point (Theorem~\ref{thm:mfg-existence}).
\vs
We conclude showing that the mean field model provides a quantitative approximation to its
finite population analogue. Among others, we show that the mean field population learning feedback is an $O((N-1)^{-1/2})$ approximate Nash equilibrium for the  population of size $N$.

\paragraph{Related literature.} There has been a substantial body of  work on mean field games with optimal stopping  in the fully
observable setting. Among others, In \cite{nutz2018mean}, Nutz studied interactions
through the number of agents who have already stopped, while Carmona,
Delarue and Lacker in \cite{carmona.delarue.lacker.2017.timing} analyzed a
bank-run game with Brownian common noise and proved existence of
mean field equilibria. The first general analytical treatment using the
coupled Hamilton--Jacobi--Bellman and Fokker--Planck system was given by
Bertucci in \cite{Bertucci}, who introduced the notion of mixed solutions.
A ``dual'' linear programming approach involving relaxed solutions was developed by Bouveret,
Dumitrescu and Tankov in \cite{bouveret.dumitrescu.tankov.20} and
Dumitrescu, Leutscher and Tankov
\cite{dumitrescu2021control,dumitrescu.leutscher.tankov.2022linear}. Further advances include a master equation
formulation by Possamai and Talbi \cite{possamai2023mean}, a mean field
Bank--El Karoui representation approach by He, Tan and Zou
\cite{he2023mean}, and related work on McKean--Vlasov optimal stopping by
Talbi, Touzi and Zhang \cite{talbi2023dynamic}, Cosso and {Perelli}
\cite{cosso2025mean}, and Cardaliaguet, Jackson and Souganidis
\cite{cardaliaguetjacksonsouganidis}.

\vs

The closest work to our formulation is  \cite{campbell2024mean}.  Their main
analysis concerns the objective $\E[c\tau+G(P_\tau)]$ and requires the
obstacle $G$ to be smooth, concave, and symmetric about $1/2$. The classical
piecewise-affine obstacle is treated separately but only under an additional
monotonicity assumption on the drift of the belief process.  Their proof of the boundary stability  directly adapts the
probabilistic approach of De Angelis and Stabile \cite{de2019lipschitz}, which necessitates 
additional structural restrictions on the data. 
Our approach, which is based on the free boundary analysis of Blanchet, Dolbeault and Monneau \cite{BDM2010} and Blanchet \cite{Blanchet2006exercise},  shows  that the first contact point
$(\tau,P_\tau)$ is almost surely regular on $\{0<\tau<T\}$; the {meaning} of the term ``regular'' is explained later in the paper. 
%(Proposition~\ref{lem:local-contact-Holder}).
This yields local $C^{0,1/2}$ regularity of the contacted boundary,
a key ingredient in establishing convergence of the optimal
stopping times.
\vs
In our model,  the initial belief being a random variable with distribution $\eta$ models the heterogeneity of prior knowledge in the large scale system. This differs from the setting of \cite{campbell2024mean}, where it is assumed that there is only common  prior belief $p_0$. 
\vs
Our objective is more general compared to \cite{campbell2024mean}. We allow for   
time and state dependent running cost and a piecewise affine cost $G$ without any additional monotonicity assumption.  %which is not necessarily symmetric about $1/2$, without any additional monotonicity assumption.  
Moreover, we establish a quantitative finite population approximation,
which is left open in \cite{campbell2024mean}.
Since the free boundary stability is
the essential continuity input in mean field games of optimal stopping, our
analysis provides a potential framework for studying a broader class of such
problems.  
%We also establish a quantitative finite-player approximation,
%which is left open in \cite{campbell2024mean}.

\paragraph{Organization of the paper.} 
%The introduction concludes  with a brief discussion about some very basic facts about 
%measures defined in bounded time intervals that we will be using throughout.
In Section~\ref{sec:model} we introduce the main assumptions,  the mean field information acquisition
model, the notion of a mean field solution and the  main result about the existence of a mean field equilibrium  as well as a number of illustrating examples. 
Section~\ref{sec:VI-boundaries} is devoted to a detailed analysis of  the underlying variational inequality
and  its stopping (free)  boundaries. 
%We construct a monotone penalization whose
%solutions converge uniformly to the value function
%(Proposition~\ref{prop:penalization}) and prove that the value function is
%concave in the belief variable (Proposition~\ref{prop:concavity}). We then
%establish the two-boundary geometry, with the exploration region lying
%between the two boundaries (Proposition~\ref{prop:two-boundary}). In Proposition~\ref{lem:local-contact-Holder} we show that the interfaces are $1/2-$H\"older continuous at the first contact points,  
%in 
%Proposition~\ref{lem:initial-piecewise} we prove that the interfaces are  are continuous at $t=0$, and, for completeness, in Proposition~\ref{prop:terminal-boundary-limits} we establish that the interfaces have a trace at $m$.  Finally,
%Theorem~\ref{thm:boundary-ae} shows  that, under the  weak convergence of the
%mean field environments, both stopping boundaries converge uniformly
%on compact subsets of a full-measure set $\mathcal R\subset(0,T)$.
In Section~\ref{sec:mfg-existence} we establish the  existence of a mean field
equilibrium. 
%For this, we first prove that the weak convergence of the mean field
%environments implies convergence in probability of the corresponding
%best response stopping-decision pairs and, thus, the continuity of the 
%best-response map $\Phi:\mathcal M\to\mathcal M$ 
%(Proposition~\ref{prop:canonical-response-continuity}). Since
%$\mathcal M $ is nonempty, compact, and convex, the
%Schauder--Tychonoff theorem yields a fixed point
%$\mu^*\in{\mathcal M}$ satisfying $\Phi(\mu^*)=\mu^*$, and hence a
%mean field equilibrium (Theorem~\ref{thm:mfg-existence}).
%\vs
Finally, in Section~\ref{sec:N-player} we formulate and analyze a finite population game and show that the mean field equilibria obtained earlier are approximate
Nash equilibria. 
\section{Preliminaries, assumptions and examples}
\label{sec:model}

We begin with some preliminaries about the weak topology of $\cM$. Then we introduce the precise assumptions and conclude with a number of examples that fit our general setting.

%We begin with a short subsection about the notation and some preliminaries. Then 
\subsection{The weak topology of $\cM$}  
%Let $(\Omega,\mathcal F,\Pp)$ be a probability space that supports an initial belief $P_0$ with law with the property that
%\beq\label{density}
%\eta \ \ \text{is   probability measure 
%on $[0,1]$ with density  denoted also by $\eta$,} 
%\eeq
%%taken to be  probability measure 
%%on $[0,1]$ with density which we denote also by $\eta$, 
%an
%unobserved binary state $\Theta\in\{0,1\}$, and a Brownian motion $W$.
%We assume % that $P_0$ has law $\eta$, a probability measure on $[0,1]$,
%that 
%\beq\label{vlakas20}
%\Pp(\Theta=1\mid P_0)=P_0 \ \ \text{and} \ \ \text{$W$ is independent of $(P_0,\Theta)$},
%\eeq
%and we 
%denote by $\mathcal M_+([0,T])$ and $\mathcal P([0,T])$ the {nonnegative Borel measures with total mass at most one}
%and  the probability measures on $[0,T]$ respectively.  
%\vs
To deal with the convergence in $\mathcal M_+([0,T])$, 
%metrize the convergence of 
we adjoin a ``cemetery point'' 
$\dagger$ to $[0,T]$,  set $\overline I=[0,T]\cup\{\dagger\}$ and introduce the distance $ d_\dagger$ on $\overline I$ given by 
\[
 d_\dagger(s,r)=|s-r|,\quad
 d_\dagger(s,\dagger)=d_\dagger(\dagger,s)=T+1 \quad \text{and} 
 \quad d_\dagger(\dagger,\dagger)=0.
\]
Let $\delta_\dagger$ be the $\delta-$function at $\dagger$. We then extend  a nonnegative  subprobability measure $\lambda$ on $[0,T]$ to a probability measure $\ol \lambda$ on {$\ol I$} by 
\[
 \bar\lambda:=\lambda +\bigl(1-\lambda([0,T])\bigr)\delta_\dagger. 
 \]
 
 %where $\delta_\dagger$ is the $\del$ function at $\dagger$. 
 In turn, the {Wasserstein} metric  $\mathcal W_1^{\mathrm{ext}}$ on $\mathcal M_+([0,T])$ is defined by 
 %\qquad
 \[
 \mathcal W_1^{\mathrm{ext}}(\lambda,\lambda')
 :=\mathcal W_1(\ol \lambda, \ol \lambda'). 
\]
%The mean field space is
%We introduce the space of sub-probability measures %The environment space, where we will find the mean field measure, is 
%\begin{equation}
% \mathcal M:=\Bigl\{(\mu^0,\mu^1):
% \mu^0,\mu^1\in\mathcal M_+([0,T]), {\mu^0+\mu^1\in \mathcal{P}([0,T])}\Bigr\},
% \label{eq:set_C}
%\end{equation}
We  equip $\cM$ with the product weak topology metrized by
\begin{equation}
 {d_{\mathcal M}}(\mu,\nu)
 :=\mathcal W_1^{\mathrm{ext}}(\mu^0,\nu^0)
  +\mathcal W_1^{\mathrm{ext}}(\mu^1,\nu^1),
 \label{eq:mean field-metric}
\nonumber
\end{equation}
and note that  it is compact and {convex}. %$\mathcal M$ is compact and convex.  
\vs
For each $t\in [0,T]$, we consider  the cumulative in time histories % In what follows we write
\[
 \mu_t^i:=\mu^i([0,t]),\]
and, for convenience, we write 
 \[\mu_{\le t}^i:=\{\mu_s^i:0\le s\le t\}  \quad \text{and} \quad \mu_{\le t}:=(\mu_{\le t}^0,\mu_{\le t}^1).
\]
In the sequel, we will refer to $\mu_{\le t}$ as collective learning input or history.

\subsection{Assumptions}
To model the gradual, rather than instantaneous, effect of the past stopping-learning decisions we need to regularize $\mu_{\le t}^0$ and $\mu_{\le t}^1$. For this, we select 
a nonnegative $\phi\in C_c^\infty((0,\infty))$  with
$\int_\R\phi(r)\dd r=1$ and introduce 
$H(r):=\int_{-\infty}^r\phi(u)\dd u$.  Then we  define, for $i=0,1$,
\begin{equation}
 F(t;\mu_{\le t})
 :=\bigl(F_0(t;\mu_{\le t}^0),F_1(t;\mu_{\le t}^1)\bigr)
\quad \text{where} \quad
 F_i(t;\mu_{\le t}^i)
 :=\int_{[0,T]}H(t-s)\,\mu^i(\dd s).
 \label{eq:causal-smoothing}
\end{equation}
Since 
$F_i(t;\mu_{\le t}^i)=\int_0^t\mu_s^i\phi(t-s)\dd s$, we 
see that 
\begin{equation}\label{takis1}
\begin{split}
& F(t;\mu_{\le t}) \ \text{depends only on stopping decisions made before
$t$, each of its components}\\
& \text{is nondecreasing in $t$, and $F$ is
$[0,1]^2$-valued and uniformly Lipschitz in $t$. } 
\end{split}
\nonumber
\end{equation}
\vs

Fix $\mu\in\mathcal M$. For the signal drift $b\bigl(t,\Theta;F(t;\mu_{\le t})\bigr)$ in \eqref{vlakas0}, with $\Delta b(t;z):=b(t,1;z)-b(t,0;z)$ for $z\in[0,1]^2$ as in \eqref{ze10}, we assume that 
\be\label{takis11}
b:[0,T]\times\{0,1\}\times[0,1]^2\to\R \ \ \text{is bounded and continuous,}
\ee

%The drift b\bigl(t,\Theta;F(t;\mu_{\le t})\bigr)
%The representative (observation) signal is modeled as 
%\begin{equation}
% \dd Y^\mu_t=b\bigl(t,\Theta;F(t;\mu_{\le t})\bigr)\dd t+\dd W_t \ \ \text{for} \ \  t\in[0,T] \ \  \textrm{ with } \ \ Y^\mu_0=0,
% \label{eq:signals}
%\end{equation}
%where $F(t;\mu_{\le t})$ is as in \eqref{eq:causal-smoothing} and 
%\be\label{takis11}
%b:[0,T]\times\{0,1\}\times[0,1]^2\to\R \ \ \text{is bounded and continuous.}
%\ee

%\vs

%Set
%\[
%\Delta b(t;z):=b(t,1;z)-b(t,0;z),
%\qquad
%\Delta b(t;\mu_{\le t})
%:=\Delta b\bigl(t;F(t;\mu_{\le t})\bigr).
%%\]
%
%Let $\big(\mathcal F_t^{Y^\mu})_{0\le t\le T}$ be the  completed filtration generated by
%$(P_0,Y^\mu_s)_{s\le t}$, and, for each $s\in [0,T]$ define the belief variable 
%\be
%P^\mu_s:=\Pp(\Theta=1\mid\mathcal F_s^{Y^\mu}). %  \ \ \text{for} \ \ s\in 
%\ee
%For
%$(t,p)\in[0,T]\times[0,1]$, we denote by $P^{t,p;\mu}$ the posterior
%process initialized at $p$ at time $t$. Classical filtering reults (see \cite{}) yield 
%that  
%%of Bain
%%and Crisan \cite[Chapter~3]{BainCrisan2009} gives
%\begin{equation}
% \dd P_s^{t,p;\mu}
% =\Delta b^\mu(s)P_s^{t,p;\mu}
%   (1-P_s^{t,p;\mu})\dd Z_s \ \ \text{for} \ \  s\in (t,T] \ \ \text{and} \ \ P_t^{t,p;\mu}=p,
% \label{eq:dynamics}
%\end{equation}
%where  $Z^\mu$ is a $\mathcal F^{Y^\mu}$-Brownian motion and $\Delta b^\mu$ is defined by \eqref{ze10}.
%%\[
% \Delta b(t;z):=\Delta b\bigl(t;F(t;\mu_{\le t})\bigr)=b(t,1;z)-b(t,0;z).
%% \qquad
%% \Delta b(t;\mu_{\le t})
%% :=\Delta b\bigl(t;F(t;\mu_{\le t})\bigr).
%\]
and %f%or the drift $ \Delta b(t;z):=\Delta b\bigl(t;F(t;\mu_{\le t})\bigr)=b(t,1;z)-b(t,0;z)$ we assume that 
\be\label{ass.b}
\begin{split}
&\text{
$[0,T]\times \mathcal M \ni (t,\mu) \to \Delta^\mu b(t)$ is Lipschitz continuous, and there exist} \\ 
&\text{
%%and the functions
%$(t,z)\mapsto\Delta b(t;z)$ and $(t,p,z)\mapsto\bar f(t,p;z)$
%are Lipschitz in their arguments.  
constants \ 
$0<\underline c_b\le\overline c_b<\infty$ \  such that} \ \ 
%$0<\underline f\le\overline f<\infty$ such that
 \underline c_b\le\Delta b^\mu(t)\le \overline c_b \ \ \text{uniformly in} \ \ [0,T]\times \mathcal M.
\end{split}
\end{equation}

\vs
The boundedness and Lipschitz conditions on $b$ and $\Delta b$ ensure that the sde in {\eqref{0}} has a unique strong solution. Since $0$ and $1$ are absorbing states, $P^{t,p;\mu}$ remains in
$[0,1]$ and is therefore a bounded $\mathbb F^\mu$-martingale.
\vs

It also follows from \eqref{ass.b} that $a^\mu$ defined in \eqref{ze1} is strictly bounded away from zero in compact subsets of $(0,1)$, since
\be\label{hob25}
a^\mu(t,p) =\frac12\big(\Delta b^\mu(t )\big)^2p^2(1-p)^2\ge \frac12 {\underline c_b^2} p^2(1-p)^2.
\ee

\vs

For the running cost function introduced  in {\eqref{eq:f-intro}}, we assume that %$f:[0,T]\times [0,1]\times [0,1]^2\to \R$,  we assume that
%\be\label{eq:f}
%f^\mu(t,p)
% :=f\bigl(t,p;F(t;\mu_{\le t})\bigr),
% \eeq
% where
$ f:[0,T]\times[0,1]\times[0,1]^2\to(0,\infty)$ 
satisfies
%\begin{equation}
% \bar f:[0,T]\times[0,1]\times[0,1]^2\to(0,\infty)\,\, \textrm{ and }\,\,
% f(t,p;\mu_{\le t})
% :=\bar f\bigl(t,p;F(t;\mu_{\le t})\bigr),
% \label{eq:f}
%\end{equation}
%with
\be \label{eq:primitive-bounds-1}
\begin{split}
& (t,p,z)\mapsto f(t,p;z) \ \text{ is Lipschitz continuous,}\\
& \text{for every $(t,z)$ the map \ 
$p\mapsto f(t,p;z)$ is concave, and } \\
&\text{there exist  constants
$0<\underline c_f\le\overline C_f<\infty$ such that, for all $(t,p,z),$}\\
& \qquad \qquad \qquad \underline c_f\leq  f(t,p;z) \leq \overline C_f. 
%\underline c_f\le f(t,p;\mu_{\le t})\le\overline c_f, \\
%&\text{and for every $(t,z)$, the map \ 
%$p\mapsto\bar f(t,p;z)$ is concave.}
\end{split}
\eeq

%The lower bound on  $f$ keeps the signal informative and information
%acquisition costly, while the concavity of $\bar f$ in the belief variable is
%used to establish concavity of the value function.
%\vs
%
%For notational convenience, we write
%\beq \label{eq:effective-forcing}
%f^\mu(t,p)=\bar f\bigl(t,p;F(t;\mu_{\le t})\bigr). % \ \text{and} \   F^\mu(t,p):=f^\mu(t,p)+A_-^\mu(t,p)G''(p),
%\eeq
%

%In what follows, we write 
%\beq\label{takis20}
%f^\mu(t,p)=\bar f\bigl(t,p;F(t;\mu_{\le t})\bigr),
%\eeq
%and 
%\begin{equation}
%       F^\mu(t,p):=f^\mu(t,p)+A_-^\mu(t,p)G''(p),
%       \label{eq:effective-forcing}
%    \end{equation}
%

%The standing conditions on the model coefficients are collected
%in Assumption~\ref{ass:coefficients} of Section~\ref{sec:VI-boundaries}.
%\vs
%As far as the stopping penalty $G:[0,1]\to \R$ 
An admissible stopping-decision pair at time $t$ consists of a
$[t,T]$-valued $\mathbb F^\mu$-stopping time $\tau$ and a $\mathcal F_\tau^\mu$-measurable decision $d\in\{0,1\}$.  Here,  $d=0$ and $d=1$ represent respectively  the choice of  the null hypothesis $\Theta=0$, and of the alternative hypothesis $\Theta=1$.
\vs
If learning terminates  at time $\tau$ with posterior belief $P^\mu_\tau$  and 
a decision $d$ is made, the resulting terminal loss $G$ is the classical cost in sequential testing problems given by  \eqref{chi4}. 
%\begin{equation}
% G(p):=\overline G(p,D(p))
%      =\min_{i=0,1}\overline G(p,i),
% \label{eq:G1}
%\end{equation}
%where, for $a_0, a_1>0$, 
%\be\label{takis12}
%\overline G(p,0)=a_0p \ \ \text{and} \ \ \overline G(p,1)=a_1(1-p), %\ \text{is convex for \ $i=0,1.$}
%\ee
%and $D$ is as in \eqref{eq:optimal_d}, where 
%\be\label{ze2}
%m=\frac{a_1}{a_0+a_1},
%\ee
%is the maximum point of $G$.
\vs
%
%$\overline G:[0,1]\times \{0,1\}\to \R$ is such that 
%\be\label{takis12}
%p\to \ol G(p,i) \ \text{is convex for \ $i=0,1.$}
%\ee
%\be\label{takis12}
%\overline G(p,0)=a_0p \ \ \text{and} \ \ \overline G(p,1)=a_1(1-p). %\ \text{is convex for \ $i=0,1.$}
%\ee
%The two loss functions in \eqref{takis12} coincide at 
%\be\label{ze2}
%m=\frac{a_1}{a_0+a_1}.
%\ee
Note that 
\be\label{Gconcave}
G \ \text{ is a concave Lipschitz continuous function with constant \ $\overline a=\max(a_0,a_1)$.}
% \ and concave}. % with a singularity at}  \ m.
\ee
For some of the arguments that follow,   it is more convenient to rewrite $G$ as 
\be\label{ze3}
 G(p)=G(m)+a_0(p-m) \ \text{if} \      p\le m\,\, \text{ and }\,\,
G(p)=G(m)-a_1(p-m) \ \text{if} \  p\ge m,
\nonumber
\ee
with $m$ given by \eqref{ze2.1}.
\vs

In order to make statements less cumbersome, we summarize  the assumptions above in 
\be\label{takis20}
\text{assume that \eqref{hob1}, \eqref{hob2}, \eqref{takis11}, \eqref{ass.b}, \eqref{eq:primitive-bounds-1} and \eqref{ze2.1} hold. }%\eqref{takis12} hold.}
% \eqref{eq:G_bar} \ and \ \eqref{ass:obstacle1}.}
\ee

\subsection{Examples within our framework}
We provide three  representative examples. % that fall within our framework.
\vs
{\bf (i)~Market liquidation.}
We consider a model for liquidity-driven information degradation in a financial  market, where 
%A population of traders learns a 
binary asset fundamentals are being learned from transaction
prices.  As trading stops and liquidation takes place, market depth deteriorates.  For some constants $ \beta_0>\beta_1>0$, we model  the informational signal to evolve according to 
\[
 \dd Y^\mu_t=\Theta\bigl(\beta_0-\beta_1(F_0(t;\mu_{\le t}^0)
                    +F_1(t;\mu_{\le t}^1))\bigr)\dd t+\dd W_t. 
%\,\, \textrm{ where }\,\, \beta_0>\beta_1>0.
\]
{As liquidation accumulates, the drift gap between the two hypotheses narrows, slowing information acquisition by the remaining traders. A constant running cost $f(t,p;\mu_{\le t})=c>0$ represents the cost of continued monitoring, so aggregate liquidation feeds back into the value of waiting for more information.}
%with aggregate exit.  We  may consider 
%having a constant opportunity running cost %$f(t,p;\mu_{\le t})=c>0$ and a terminal
%and a loss $G$ associated with  trading decisions.
%\end{example}
\vs
{\bf (ii)~Selective disclosure.}
%\begin{example}[Review-driven learning with selective disclosure]
We consider a review-driven learning with selective disclosure. Consumers decide whether a product is of high or low quality while reviews from earlier buyers affect both the signal and the learning running cost. %cost of waiting.
For example, we may assume that, for some constants $\alpha,\lambda,\kappa>0$, the  signal   $Y^\mu$ evolves according to 
\[
 \dd Y^\mu_t=\alpha\bigl(\Theta+\lambda F_0(t;\mu_{\le t}^0)
                   -\kappa F_1(t;\mu_{\le t}^1)\bigr)\dd t+\dd W_t.
%\,\, \textrm{ with }\,\,\alpha,\lambda,\kappa>0.
\]
The aggregate reviews shift the observed signal level, while
$\Delta b=\alpha$.  A running cost such as
$$f(t,p;\mu_{\le t})=c_0+c_1(F_0(t;\mu_{\le t}^0)
 +F_1(t;\mu_{\le t}^1)),$$ where $c_0>0$ and $c_1\ge0$, captures
declining attention or availability.
%\end{example}
\vs

%\begin{example}[Informational crowding in platform experimentation and A/B testing]
{\bf {(iii)}~Information crowding.} In  experimentation platforms with A/B testing, many product teams experiment in parallel. However,  widespread adoption reduces
experimental diversity.  We model, for some constants $\beta,\gamma>0$,  the  signal to evolve according to  
\[
 \dd Y^\mu_t=\Theta\beta\exp\big (-\gamma(F_0(t;\mu_{\le t}^0)
                         +F_1(t;\mu_{\le t}^1))\big)\dd t+\dd W_t. \,\,
%\textrm{ for }\,\, \beta,\gamma>0.
\]

We then see  that the drift decreases as more teams conclude their experiments.
  A cost
$$f(t,p;\mu_{\le t})=c_0+c_1t,$$
with $c_0>0$ and $c_1\ge0$, represents
the opportunity cost of delaying deployment.  Among others, such feedback may generate
early stopping and late experimentation equilibria.
%\end{example}

%For later use, we introduce the PDE shorthand
%\begin{equation}
% \gammA_-^\mu(t):=\Delta b(t;\mu_{\le t}),
% \qquad
% f^\mu(t,p):=f(t,p;\mu_{\le t}).
% \label{eq:primitive-coefficients}
%\end{equation}
%
\section{The variational inequality {and} the free boundaries}
\label{sec:VI-boundaries}

This section is devoted to the analysis of the variational inequality, %\eqref{eq:HJB1-0},  
the properties of its solutions and their connection with the value function, as well as some important properties of  the free boundaries. 
%  and their stability properties with respect to the weak convergence of measures in $\cM$. 

\subsection{Stability  with respect to the collective input}

In view of the particular dependence of the drift of the signal and the running cost on the collective input, 
we have the following result.

\begin{proposition}%[Continuity of the coefficients]
\label{lem:coefficient-stability}
If $\mu^n\weakto\mu$ in {$\mathcal M$}, then
\begin{equation}
 \|{\Delta b^{\mu^n}-\Delta b^\mu}\|_{L^\infty(0,T)}
% \|\Delta b(s;\mu^n_{\le s})-\Delta b(s;\mu_{\le s})\|_{L^\infty(0,T)}
 +\|f^{\mu^n}-f^\mu\|_{L^\infty([0,T]\times[0,1])}\longrightarrow0.
 \label{eq:data-stability}
\end{equation}
\end{proposition}

\begin{proof}
For every $t$, the function $s\mapsto H(t-s)$ on $[0,T]$, extended by
the value $0$ at $\dagger$, is Lipschitz on $(\overline I,d_\dagger)$
with constant at most
$\max\{\|\phi\|_\infty,(T+1)^{-1}\}$.  
\vs
It then follows from the
Kantorovich--Rubinstein duality that
\[
 \begin{split}
 \sup_{t\in[0,T]}
 |F\Big(t; \big((\mu^n)_{\le t}^0, (\mu^n)_{\le t}^1\big)\Big )&-F\Big(t;\big(\mu_{\le t}^0, \mu_{\le t}^1\big)\Big)|
\\
& \le \max\{\|\phi\|_\infty,(T+1)^{-1}\}
       \big[\mathcal W_1^{\mathrm{ext}}(\mu^{n,0},\mu^0)+ \mathcal W_1^{\mathrm{ext}}(\mu^{n,1},\mu^1)\big]
 \longrightarrow0.
\end{split}
\]
The Lipschitz continuity of $\Delta b$ and $\bar f$ in the aggregate
variable proves \eqref{eq:data-stability}.

\end{proof}

\subsection{Viscosity solutions}

%\begin{lemma}[Spatial and terminal boundary values]
%\label{lem:boundary-values}
%For every $t\in[0,T]$,
%\begin{equation}
%    V^\mu(t,0)=G(0),\qquad V^\mu(t,1)=G(1),
%    \label{eq:spatial-boundary-values}
%\end{equation}
%and, for every $p\in[0,1]$,
%\begin{equation}
%    V^\mu(T,p)=G(p).
%    \label{eq:terminal-value}
%\end{equation}
%\end{lemma}
%
%\begin{proof}
%If $p=0$, then $P_s^{t,0;\mu}=0$ for all $s\ge t$.  Consequently, for
%every $\tau\in\mathcal T_{t,T}$,
%\[
% \E\left[\int_t^\tau f^\mu(s,0)\dd s+G(0)\right]\ge G(0),
%\]
%and equality is attained by $\tau=t$.  The proof at $p=1$ is identical.
%When $t=T$, the only admissible stopping time is $T$, which proves
%\eqref{eq:terminal-value}.
%\end{proof}

We begin with the definition of the viscosity solutions to the terminal value problem \eqref{takis100}.

\begin{definition}\label{def:viscosity}
(i)~A bounded {continuous} function $u:[0,T]\times[0,1]\mapsto \mathbb{R}$ is  a viscosity subsolution to {\eqref{takis100}},     if  $u(T,\cdot)\le G(\cdot)$ on $[0,1]$ and, for any test function $\phi\in C^{1,2}([0,T]\times[0,1])$ and any $(t_0,p_0)\in (0,T)\times[0,1]$, such that $u-\phi$ attains a local maximum at $(t_0,p_0)$ with $u(t_0,p_0)=\phi(t_0,p_0)$, 
    \begin{eqnarray}\label{eq:viscosity-sub}
   \max \Big\{- \partial_t \phi(t_0,p_0) - a^\mu(t_0,p_0;\mu_{\leq t})\pt^2_p \phi(t_0,p_0) - f(t_0,p_0;\mu_{\leq t}),  u(t_0,p_0)- G(p_0)\Big\}\leq 0.
    \nonumber
    \end{eqnarray}
    %\item (Terminal condition) For all $p\in[0,1]$, $u(T,p)\leq G(p)$.
 %   \item (Boundary condition) For all $t\in[0,T]$, $u(t,0)\leq G(0)$ and  $u(t,1)\leq G(1)$.

(ii)~A bounded {continuous} function $w:[0,T]{\times}[0,1]\rightarrow \mathbb{R}$ is a viscosity supersolution to {\eqref{takis100}}, if $w(T,\cdot)\ge G(\cdot)$ on $[0,1]$ and   for any test function $\phi\in C^{1,2}([0,T]\times[0,1])$ and 
any $(t_0,p_0)\in (0,T)\times(0,1)$ such that $w-\phi$ attains a local minimum  with $w(t_0,p_0)=\phi(t_0,p_0)$, 
    \begin{eqnarray}\label{eq:viscosity-super}
\max \Big\{- \partial_t \phi(t_0,p_0) - a^\mu(t_0,p_0;\mu_{\leq t})\partial^2_{p} \phi(t_0,p_0) -f(t_0,p_0;\mu_{\leq t}),    w(t_0,p_0)-G(p_0)\Big\}\geq 0.
    \nonumber
    \end{eqnarray}
(iii) A bounded continuous function $v$ is a viscosity solution to  {\eqref{takis100}},  if it is both a subsolution and a supersolution.
\end{definition}

To simplify the exposition, we omit from now on the term ``viscosity'' with the understanding that all (subsolutions and supersolutions) are in the viscosity sense. 

\vs 
We note that, in view of the explicit form of $a^\mu$ at $p=0$ and $p=1$, no boundary conditions are needed, meaning that  subsolutions and  supersolutions to {\eqref{takis100}} are automatically subsolutions and supersolutions  at $p=0$ and $p=1$.  This is the result of the next lemma.

\begin{lemma}%[No spatial boundary condition]
\label{lem:no-spatial-boundary-condition}
Fix  $\mu\in\mathcal M$ and  let
$u,w\in C([0,T]\times[0,1])$ be, respectively, a 
subsolution and a  supersolution to {\eqref{takis100}} in
$(0,T)\times(0,1)$. Then $u$ is a 
subsolution  and $w$ is  a 
supersolution to {\eqref{takis100}} in
$[0,T)\times[0,1]$.
\end{lemma}

\begin{proof}
Since both claims are proved similarly, we only present the details
for the subsolution case.
\vs
We also note  that we only need to check what happens if, for a test function $\phi\in C^{1,2}([0,T]\times[0,1])$,  $u-\phi$ attains a maximum at $(t_0,p_0)$ for $t_0\in [0,T)$ and $p_0\in\{0,1\}$ and $u(t_0,p_0)=\phi(t_0,p_0)$.
% or $(t_0,1)$ for $t_0\in [0,T)$ and $u(t_0,1)=\phi(t_0,1)$. 
\vs
We first consider  the former case and note that, without any loss of generality, we may assume that the maximum is strict. If not,  we replace $\phi$ by $\phi(t,p)+K\bigl((t-t_0)^2+p^2\bigr)^2$
for a sufficiently large $K>0$.
\vs
%that, for some test function
%$\phi\in C^{1,2}([0,T]\times[0,1])$, $u-\phi$ attains a local maximum
%at $(t_0,0)$ with $t_0\in[0,T)$ and
%$u(t_0,0)=\phi(t_0,0)$, and note that we may assume replacing
%$\phi$ by $\phi(t,p)+K\bigl((t-t_0)^2+p^2\bigr)^2$
%for a sufficiently large $K>0$ if necessary,  that the maximum is global and strict.
%% (by replacing
%%$\phi$ by $\phi(t,p)+K\bigl((t-t_0)^2+p^2\bigr)^2$
%%for a sufficiently large $K>0$ if necessary).
%\vs
Then, the map
\[
 (t,p)\longmapsto
 u(t,p)-\phi(t,p)
 -\varepsilon\left(\frac1p+\frac1t\right)
\]
attains a maximum at some
$(t_\varepsilon,p_\varepsilon)\in(0,T)\times(0,1)$.
\vs
It then follows from the {continuity} of $u$ and the strictness of the maximum above  that, 
as $\varepsilon\to0$,
\begin{equation}
 (t_\varepsilon,p_\varepsilon)\to(t_0,0)\,\,
\textrm{ and }\,\,
 \varepsilon\left(
 \frac1{p_\varepsilon}+\frac1{t_\varepsilon}
 \right)\to0.
 \label{aus1}
\end{equation}

Since $u$ is a subsolution in $(0,T)\times(0,1)$, it follows that
\[
\max\Bigg\{
 -\partial_t\phi(t_\ep,p_\ep)+\frac{\ep}{t_\ep^2}
 -a^\mu(t_\ep,p_\ep;\mu_{\leq t_\ep})
 \left(\pt_p^2\phi(t_\ep,p_\ep)
       +\frac{2\ep}{p_\ep^3}\right)
 -f(t_\ep,p_\ep;\mu_{\leq t_\ep}),
 \ u(t_\ep,p_\ep)-G(p_\ep)
\Bigg\}\leq0.
\]
Using \eqref{ze1}, we find that
\[
 a^\mu(t_\ep,p_\ep;\mu_{\leq t_\ep})
 \frac{2\ep}{p_\ep^3}
 =
 \frac{\ep}{p_\ep}
 \Bigl(\Delta b(t_\ep;
                 \mu_{\leq t_\ep})\Bigr)^2
 (1-p_\ep)^2,
\]
and, hence, 
\begin{equation}
 \lim_{\ep\to0}
 a^\mu(t_\ep,p_\ep;\mu_{\leq t_\ep})
 \frac{2\ep}{p_\ep^3}=0.
 \label{aus2}
\end{equation}

Moreover, by the continuity of the coefficients, we have 
\begin{equation}
\begin{aligned}
 &\lim_{\ep\to0}
 a^\mu(t_\ep,p_\ep;\mu_{\leq t_\ep})
 \pt_p^2\phi(t_\ep,p_\ep)
 =
 a^\mu (t_0,0;\mu_{\leq t_0})\pt_p^2\phi(t_0,0),\\
 &\lim_{\ep\to0}
 f(t_\ep,p_\ep;\mu_{\leq t_\ep})
 =
 f(t_0,0;\mu_{\leq t_0}).
\end{aligned}
\label{aus3}
\end{equation}

Using \eqref{aus1}, \eqref{aus2} and \eqref{aus3}, and omitting  the
nonnegative term $\eps/t_\eps^2$, we obtain the sub-solution inequality
at $(t_0,0)$. %The other cases are analogous
\vs

If the maximum occurs at $p_0=1$, we repeat the same argument the only difference being that we use $\frac{1}{1-p}$ in place of $\frac{1}{p}$.

\end{proof}

The following result is well known  and  we refer to  Crandall, Ishii and Lions \cite{CIL} for its proof.

\begin{theorem}\label{uniq}
Let $u$ and $w$ be respectively {a bounded continuous subsolution and a bounded continuous supersolution} to {\eqref{takis100}}, and assume that $u(T,\cdot)\leq w(T,\cdot)$ on $[0,1]$. Then, $u\le w$ on $[0,T]\times [0,1]$. In particular, 
{\eqref{takis100}} has a unique continuous viscosity solution.
\end{theorem}

Next, we argue that, for each $\mu \in \cM$, the value function $V^\mu$ is the unique continuous viscosity solution to {\eqref{takis100}}. 
A key fact towards establishing this claim is that the value function satisfies the dynamic programming principle, which states (see Fleming and Soner \cite[Chapters~V and VIII]{FlemingSoner2006}) that,
for every $\theta\in\mathcal T_{t,T}$,
\beq\label{dp}
 V^\mu(t,p)=\inf_{\tau\in\mathcal T_{t,T}}\E\bigg[
 \int_t^{\tau\wedge\theta}f^\mu(s,P_s^{t,p;\mu})\dd s + V^\mu\big(\tau \wedge \theta , P_{\tau\wedge \theta}^{t,p;\mu}\big)
  \bigg].
% 
% 
% +\1_{\{\tau<\theta\}}G(P_\tau^{t,p;\mu})+\1_{\{\tau\ge\theta\}}
% V^\mu(\theta,P_\theta^{t,p;\mu})\bigg].
\eeq

We continue with the space-time continuity property of $V^\mu$.
% first consider the continuity of the value function and then we show that it is indeed a viscosity solution to the variational inequality. 

\begin{proposition}%[Continuity and viscosity characterization]
\label{prop:viscosity-characterization}
Assume \eqref{takis20} and fix $\mu \in \cM$. For each $t\in [0,T]$,  $V^\mu(t,0)=G(0)$ and $V^\mu(t,1)=G(1)$ and, for each $p\in [0,1]$, $V^\mu(T,p)=G(p)$. Moreover, $V^\mu$ is Lipschitz  continuous in space and $1/2-$H\"older continuous in time on $[0,T]\times [0,1]$.
% and is the
%unique continuous viscosity solution of
%\begin{equation}
% \max\Big\{-v_t-A_-^\mu v_{pp}-f^\mu,\,v-G\Big\}=0
% \quad\text{in }[0,T)\times[0,1], \,\,\textrm{ with } \,\,v(T,p)=G(p)\,\, \textrm{ for } \,\,p\in[0,1]. 
% \label{eq:VI}
%\end{equation}
%
\end{proposition}

\begin{proof}
If $p=0$, then $P_s^{t,0;\mu}=0$ for all $0\leq t \leq s\le T$.  Consequently, for
every $\tau\in\mathcal T_{t,T}$,
\[
 \E\left[\int_t^{\tau\wedge T} f^\mu(s,0)\dd s+G(0)\right]\ge G(0),
\]
and equality is attained by $\tau=t$.  The proof at $p=1$ is similar.
When $t=T$, the only admissible stopping time is $T$, and, hence, $V^\mu(T,p)=G(p)$ for all $p\in [0,1]$. 
%which proves \eqref{eq:terminal-value}.
\vs
To prove the Lipschitz continuity in space, 
take $t\in[0,T]$ and $p,q\in[0,1]$, and recall that %and construct
$P^{t,p;\mu}$ and $P^{t,q;\mu}$ solve  {\eqref{0}} with the same Brownian motion.
Since, in view of \eqref{ze1} and \eqref{hob25},  the diffusion coefficient is uniformly Lipschitz in the state variable away from $\{0,1\}$, 
the Burkholder--Davis--Grundy and Gr\"onwall inequalities give
\begin{equation}
\E\left[
\sup_{t\le r\le T}
\left|P_r^{t,p;\mu}-P_r^{t,q;\mu}\right|^2
\right]
\le C|p-q|^2.
\label{eq:sde-space-estimate}
\nonumber
\end{equation}
Hence, for every $\tau\in\mathcal T^\mu_{t,T}$,
\begin{align}
&\left|
\E\left[
\int_t^\tau f^\mu(r,P_r^{t,p;\mu})\dd r
+G(P_\tau^{t,p;\mu})
\right]
-
\E\left[
\int_t^\tau f^\mu(r,P_r^{t,q;\mu})\dd r
+G(P_\tau^{t,q;\mu})
\right]
\right| \notag\\
&\qquad\le
L_f\int_t^T
\E\left|P_r^{t,p;\mu}-P_r^{t,q;\mu}\right|\dd r
+
\overline a \E\left|
P_\tau^{t,p;\mu}-P_\tau^{t,q;\mu}
\right| \le C|p-q|, \notag
%&\qquad\le C|p-q|.
\label{eq:cost-space-estimate}
\end{align}
where $\ol a$  is as in \eqref{Gconcave}.
\vs
Taking infima over $\tau\in\mathcal T^\mu_{t,T}$ in the above estimates gives
\begin{equation}
|V^\mu(t,p)-V^\mu(t,q)|\le C|p-q|.
\label{eq:value-space-continuity}
\end{equation}

For the continuity in time, we fix $0\le t<s\le T$, use \eqref{dp} with  $\theta=s$ and choose an $\ep-$approximate minimizing stopping time $\tau^\ep$  for $V^\mu(t,p)$. Using that by, its definition, $V^\mu(s,p)\leq G(p)$,  the positivity of $f^\mu$, the Lipschitz continuity in space of $G$ and $V^\mu$ and the fact that, for some $c>0$, 
\[\E\big[|P_\rho^{t,p;\mu} -p|\big]\le c \sqrt{|\rho-t|},\]
we find
\[
\begin{split}
& V(t,p) -V(s,p) \ge \E\bigg[
 \int_t^{\tau^\ep\wedge s}f^\mu(\rho,P_\rho^{t,p;\mu})\dd \rho +\1_{\{\tau^\ep<s\}}G(P_{\tau^\ep}^{t,p;\mu})\\[1.5mm]
 &\hskip1.75in +\1_{\{\tau^\ep\ge s\}}
 V^\mu(\theta,P_s^{t,p;\mu})-
V^\mu(s, p)\bigg] - \ep \\[1.5mm]
 &\qquad \qquad \ge  \E\big[\1_{\{\tau^\ep<s\}}\big(G(P_{\tau^\ep}^{t,p;\mu})-G(p)\big)\big] +
\E\big[\1_{\{\tau^\ep\ge\theta\}}\big(V^\mu(s,P_s^{t,p;\mu})-V^\mu(s, p)\big)\big] -\ep\\[1.5mm]
& \qquad \qquad \ge -\ol a \E\big[\1_{\{\tau^\ep<s\}}| P_{\tau^\ep}^{t,p;\mu}-p|\big] -C\E\big[|P_s^{t,p;\mu}-p|\big] +\ep \\[2.5mm]
&  \qquad \qquad \ge -c(\ol a +C)\sqrt{|s-t|} -\ep. 
 \end{split}
 \]
A similar argument,   choosing now an $\ep-$minimizing stopping time for $V^\mu(s,p)$ leads, for some $C'>0$, to 
\[V(s,p)-V(t,p)\ge -C'\big (L_f |s-t| + (\ol a+C) |s-t|^{1/2}\big) -\ep.\] 
It then follows, after letting $\ep\to 0$, that  
\beq\label{eq:value-time-continuity} 
|V(s,p)-V(t,p)| \lesssim \big(|s-t|+\sqrt{|s-t|}\big).
\eeq

\end{proof}

We explore next the connection between  $V^\mu$ and the solution to  {\eqref{takis100}}.

\begin{proposition}\label{viscs}
Assume \eqref{takis20}. For any $\mu \in \cM$, $V^\mu$ is the unique continuous viscosity solution to {\eqref{takis100}}. \end{proposition}

\begin{proof}
In view of Theorem~\ref{uniq} and Proposition~\ref{prop:viscosity-characterization}, we only need to check that $V^\mu$ is a viscosity solution. 
\vs
It follows from \eqref{dp} %The dynamic programming principle states
that, for every $\theta\in\mathcal T_{t,T}$,
\begin{align*}
 V^\mu(t,p)=\inf_{\tau\in\mathcal T_{t,T}}\E\bigg[&
 \int_t^{\tau\wedge\theta}f^\mu(s,P_s^{t,p;\mu})\dd s
 +\1_{\{\tau<\theta\}}G(P_\tau^{t,p;\mu})+\1_{\{\tau\ge\theta\}}
 V^\mu(\theta,P_\theta^{t,p;\mu})\bigg].
\end{align*}
To obtain the viscosity inequalities at $(t,p)$, we introduce, for $h,r>0$ small, the stopping time 
\[
 \theta_h:=(t+h)\wedge
 \inf\{s\ge t:|P_s^{t,p;\mu}-p|\ge r\}.
\]
%with  $h,r>0$  small.  
If a smooth function $\phi$ touches $V^\mu$ from above at $(t,p)$,
using the dynamic programming principle  and It\^o's formula gives
$$-\pt_t\phi-a^\mu\pt^2_p\phi -f^\mu\le0.$$

Since immediate stopping gives
\(V^\mu-G\le0\), we obtain
\[
 \max\{-\pt_t\phi-a^\mu\pt^2_p\phi -f^\mu,\,
        V^\mu-G\}(t,p)\le0.
\]
If $\phi$ touches $V^\mu$ from below at $(t,p)$, then
\[
\max\{-\pt_t\phi-a^\mu\pt^2_p\phi -f^\mu,\,
       V^\mu-G\}(t,p)
\ge
\begin{cases}
-\pt_t\phi-a^\mu\pt^2_p\phi -f^\mu\ge0 \ \ & \text{if} \ \ 
     V^\mu(t,p)<G(p),\\[1mm]
0 \ \ & \text{if} \ \ 
     V^\mu(t,p)=G(p).
\end{cases}
\]

\end{proof}
%Thus $V^\mu$ solves \eqref{eq:VI}. At \(p=i\in\{0,1\}\), test
%functions touch relative to \([0,T)\times[0,1]\), and
%\[
% P_s^{t,i;\mu}=i,\,\,\,\,
% V^\mu(t,i)=G(i),\,\,\,\,
% A_-^\mu(t,i)=0,\,\,\textrm{ for }\,\, s\ge t.
%\]
%Together with Lemma \ref{lem:no-spatial-boundary-condition}, this proves that \(V^\mu\) solves
%\eqref{eq:VI}.
%
%
%The comparison principle
%\cite[Chapters~V and VIII]{FlemingSoner2006} applies because
%\[
% A_-^\mu,f^\mu\in C([0,T]\times[0,1]),\,\,
% G\in C([0,1]),\,\, A_-^\mu\ge0,
%\]
%and
%\[
% \bigl|\gammA_-^\mu(t)p(1-p)-\gammA_-^\mu(t)q(1-q)\bigr|
% \le\overline\gamma|p-q|.
%\]
%Thus, if \(u\) is a bounded upper semicontinuous viscosity subsolution
%of \eqref{eq:VI} and \(w\) is a bounded lower semicontinuous viscosity
%super-solution, then
%\[
% u(T,p)\le w(T,p)\quad\text{for every }p\in[0,1]
% \quad\Longrightarrow\quad
% u(t,p)\le w(t,p)\quad\text{for every }(t,p)\in[0,T]\times[0,1].
%\]
%Consequently, if \(u\) and \(w\) are continuous viscosity solutions
%with terminal value \(G\), then
%\[
% u\le w,\qquad w\le u,\,\,\text{and hence}\,\, u=w.
%\]
%
%
%\end{proof}

The next result is an immediate consequence of the stability of viscosity solutions and the continuity of the coefficient of {\eqref{takis100}} with respect to the weak convergence in $\cM$.

\begin{proposition}%[Uniform stability of the value functions]
\label{prop:value-stability}
If $\mu^n\weakto\mu$ in $\cM$, then
\begin{equation}
 \|V^{\mu^n}-V^\mu\|_{L^\infty([0,T]\times[0,1])}\longrightarrow0.
 \label{eq:value-stability}
\nonumber
\end{equation}
\end{proposition}

\begin{proof}
We first observe that the estimates \eqref{eq:value-space-continuity} and
\eqref{eq:value-time-continuity} are uniform in $\mu$.  Hence,  the sequence 
$(V^{\mu^n})_n$ is uniformly bounded and equicontinuous on
$[0,T]\times[0,1]$.  Let $(V^{\mu_{n_k}})_k$ be any subsequence.
By the Arzel\`a--Ascoli theorem, it has a further subsequence converging
uniformly to a  continuous function $\overline V$.
\vs
Since, as claimed in Proposition~\ref{lem:coefficient-stability}, %$ \eqref{eq:def-coefficients} and \eqref{eq:data-stability},
\[
\|a^{\mu_{n_k}}-a^\mu\|_\infty \ \ \text{and} \ \ 
 \|f^{\mu_{n_k}}-f^\mu\|_\infty\longrightarrow0,
\]
the stability of viscosity solutions implies that $\overline V$ solves
{\eqref{takis100}}  with the same terminal data.  
\vs
Then, in view of the uniqueness in
Proposition~{\ref{viscs}}, we have $\overline V=V^\mu$.
Thus, every subsequence has a further subsequence converging uniformly to
$V^\mu$, and we easily conclude.
%which proves \eqref{eq:value-stability}.

\end{proof}

\subsection{Concavity of the value function in the space variable}

To establish the concavity of the value function in the belief variable, we consider a standard  penalization of the variational inequality. The reason is that it provides an equation instead of a variational  inequality which can be, in turn,  differentiated to obtain estimates. 
\vs
Let 
$\beta\in C^\infty(\R)$ be %convex and nondecreasing, with
\begin{equation}
\begin{split}
&\text{convex and nondecreasing, \ with}\\
& \beta(r)=0 \quad\text{for }r\le0,
 \quad
 \beta(r)>0 \quad\text{for }r>0,
 \quad
 \lim_{r\to\infty}\beta(r)=\infty \quad \text{and} \quad  \beta(1)\ge\overline C_f, 
 \end{split}
 \label{eq:penalty-beta}
\end{equation}
and, for $\eps>0$, set
\begin{equation}
 \beta_\eps(r):=\beta\left(\frac{r}{\eps}\right).
 \label{eq:penalty-function}
\end{equation}

For fixed $\mu\in\cM$, we then consider the auxiliary problem
\begin{equation}
 \begin{cases}
 -\pt_tV^{\mu,\eps}
 -a^\mu \pt^2_p V^{\mu,\eps}
 -f^\mu
 +\beta_\eps\bigl(V^{\mu,\eps}-G\bigr)=0 \ 
 & \text{in} \ \ [0,T)\times[0,1],\\[1.2mm]
 V^{\mu,\eps}(T,\cdot)=G \
 & \text{on} \ \ [0,1].
 \end{cases}
 \label{eq:penalized-equation}
\end{equation}
At $p=0,1$, \eqref{eq:penalized-equation} is understood in the
viscosity sense relative to the closed state space  without imposing any boundary
conditions.
\vs
The following result is well known, so we only provide a sketch of its proof.
\begin{proposition}%[Penalization limit]
\label{prop:penalization}
If $0<\eps'<\eps$, then
\begin{equation}
 V^\mu\le V^{\mu,\eps'}\le V^{\mu,\eps}
 \quad\text{on }\quad [0,T]\times[0,1],
 \label{eq:penalty-order}
\nonumber
\end{equation}
and %as $\ep\to 0$, 
\begin{equation}
 \|V^{\mu,\eps}-V^\mu\|_{L^\infty([0,T]\times[0,1])}
 {\le\eps}%\longrightarrow0.
 %\qquad\text{as }\eps\to 0.
 \label{eq:penalty-convergence}
\nonumber
\end{equation}
\end{proposition}

\begin{proof}
Since  $\beta_\eps$ is nondecreasing,  the standard viscosity existence theorem gives a bounded
continuous solution of \eqref{eq:penalized-equation}, while  the comparison
principle yields uniqueness; 
 see
\cite[Theorem~4.1 and Remark~4.5]{CIL}.
\vs
Since $ \beta_{\eps'}(r)\ge\beta_\eps(r)$ if  $0<\eps'<\eps$, 
%then,  for every $r\in\R$, 
%\[
% \beta_{\eps'}(r)\ge\beta_\eps(r), 
%\]
it follows that 
%Indeed, both sides vanish when $r\le0$, while for $r>0$ this follows from
%$r/\eps'>r/\eps$ and the monotonicity of $\beta$.  Thus, at any lower
%test function for $V^{\mu,\eps}$, the left-hand side of the
%$\eps'$-equation is the corresponding expression for the
%$\eps$-equation plus
%\[
% \beta_{\eps'}(V^{\mu,\eps}-G)
% -\beta_\eps(V^{\mu,\eps}-G)\ge0.
%\]
$V^{\mu,\eps}$ is a viscosity supersolution of the
$\eps'$-equation, and, thus, 
\[V^{\mu,\eps'}\le V^{\mu,\eps}.\]

We next compare $V^{\mu,\eps}$ with $V^\mu$ noting that, in view of the fact that \(V^\mu\le G\), \eqref{eq:penalty-beta} gives $\beta_\eps\bigl(V^\mu-G\bigr)=0$. Using  that 
\[-\pt_t V^\mu-a^\mu \pt_p^2 V^\mu-f^\mu\leq 0 \ \text{in} \ [0,T)\times [0,1],\]
we deduce that $V^\mu$ is a sub-solution to {\eqref{eq:penalized-equation}} and, thus, for each $\ep>0$, 
\[ V^\mu \le V^{\mu,\ep} \ \text{on} \ [0,T]\times [0,1].\]
%Hence, whenever \(\phi\) touches \(V^\mu\) from above at \((t,p)\),
%\[
% -\phi_t-A_-^\mu\phi_{pp}-f^\mu
% +\beta_\eps\bigl(V^\mu-G\bigr)
% =-\phi_t-A_-^\mu\phi_{pp}-f^\mu\le0.
%\]
%Thus \(V^\mu\) is a viscosity subsolution of
%\eqref{eq:penalized-equation}.
Next, we claim that $U=V^\mu +\ep$ is a supersolution  to \eqref{eq:penalized-equation}, a fact, which in combnation 
with the previous bound,  yields \eqref{eq:penalty-convergence}. 
\vs
Assume that, for a smooth $\phi$,  $U-\phi$ attains a local minimum at some $(t_0,p_0) \in [0,T)\times (0,1)$. It follows that, at  $(t_0,p_0)$,
\beq\label{japan1}
-\pt_t\phi -a^\mu \pt^2_p\phi - f^\mu + \beta_\ep(U-G)\geq 0,
\eeq
and, thus, the claim.
\vs

Indeed, if $V^\mu(t_0,p_0)<G(p_0)$, then, since $V^\mu-\phi$ attains a minimum at $(t_0,p_0)$, the supersolution property of $V^\mu$ gives that, at $(t_0,p_0)$, 
\[-\pt_t\phi -a^\mu \pt^2_p\phi - f^\mu \geq 0, \]
and, in view of the nonnegativity of $\beta_\ep$, \eqref{japan1} follows.
\vs 
If $V^\mu(t_0,p_0)=G(p_0)$, then, we have that $(t_0,p_0)$ is a minimum of $G-\phi$. It follows that $p_0\neq m$ and, in addition, 
\[\pt_t\phi (t_0,p_0)-a^\mu \pt^2_p\phi (t_0,p_0)\geq 0.\]
Moreover, in view of \eqref{eq:penalty-beta}, at $(t_0,p_0)$, 
\[-f^\mu +\beta_\ep(U-G)=-f^\mu +\beta_\ep(V^\mu-G)\geq -\ol C_f +\beta_\ep(\ep)=-\ol C_f +\beta (1)\geq 0.\]
Combining the last two inequalities, yields \eqref{japan1}.
\vs
The proof is now complete.

\end{proof}

Next, we state and prove the claim about the spatial concavity of the value function.  
\begin{proposition}%[Concavity]
\label{prop:concavity}
Assume \eqref{takis20}.  
For every $\mu\in\cM$ and each $t\in[0,T]$, the map
$p\mapsto V^\mu(t,p)$ is concave on $[0,1]$.
\end{proposition}

\begin{proof}
We prove the result by showing the claim for the solution to the auxiliary problem  \eqref{eq:penalized-equation}. To deal with the additional difficulty that $G$ is not and $f^\mu$ may not be in $C^2$, we approximate them  by  smooth concave functions $G^\delta$ and $f^{\mu,\delta}$. The stability of the solutions to the penalized problem yield that the solution $V^{\mu,\delta}$ to  \eqref{eq:penalized-equation} with obstacle $G^\delta$ and running cost $f^{\mu,\delta}$ converges uniformly to the solution of the problem with obstacle $G$ and running cost $f^{\mu}$. Hence, we only show that solutions to  the regularized \eqref{eq:penalized-equation}, which for simplicity we also denote by $V^{\mu}$, are concave. 
\vs

%In what follows to simplify the notation, we omit any dependence on the approximation parameters $\delta$ and $\ep$ and the measure $\mu$ and simply write $u$ for $V^{\mu, \ep, \delta}$/. 
%
%And for this we consider a regularized version with Fix $\mu\in\mathcal C$ and $\eps>0$. Under either alternative in
%Assumption~\ref{ass:obstacle}, extend $G$ concavely outside $[0,1]$
%and extend $f^\mu(t,\cdot)$ concavely outside $[0,1]$, and approximate
%both by convolution. This gives smooth functions $G_n$ and $f_n^\mu$ such that
%\[
% G_n''\le0,\,\, \pt^2_pf_n^\mu\le0,\,\, \textrm{ and }
% \|G_n-G\|_\infty+\|f_n^\mu-f^\mu\|_\infty\rightarrow 0 \textrm{ as } n \to \infty.
%\]
%Let $V_n^{\mu,\eps}$ solve \eqref{eq:penalized-equation} with
%$G_n$ and $f_n^\mu$ in place of $G$ and $f^\mu$, and set
%\[
% d_n(t):=\|G_n-G\|_\infty
% +(T-t)\|f_n^\mu-f^\mu\|_\infty.
%\]
%Since $\beta_\eps$ is nondecreasing,
%$V_n^{\mu,\eps}-d_n$ and $V_n^{\mu,\eps}+d_n$ are, respectively, a
%subsolution and a super-solution of \eqref{eq:penalized-equation}.
%Comparison principle therefore gives
%\begin{eqnarray}\label{eq:Vto0}
%     \|V_n^{\mu,\eps}-V^{\mu,\eps}\|_\infty
% \le d_n(0)\to 0 \textrm{ as } n \to \infty.
%\end{eqnarray}
%It is thus enough to prove that $V_n^{\mu,\eps}(t,\cdot)$ is concave.

To ease the notation, we write $w:=\pt^2_p V^{\mu}$.  After  differentiating the penalized equation
twice with respect to $p$ and using the convexity of $\beta$ {and the concavity of} $f^\mu$ and $G$, we get  
\[
\begin{aligned}
&-\pt_t w-a^\mu \pt^2_p w-2\pt_p a^\mu w_p
 +\Bigl(-\pt^2_p a^\mu+
 \beta_\eps'(V^{\mu,\eps}-G)\Bigr)w \\[1.2mm]
&\qquad =
 \pt^2_pf^\mu
 +\beta_\eps'({V^{\mu,\eps}}-G)\pt^2_p G
 -\beta_\eps''({V^{\mu,\eps}}-G)
   (\pt_pV^{\mu,\eps}-\pt_pG)^2
 \le0.
\end{aligned}
\]

Since $w(T,y)=\pt^2_pG(y)\leq 0$, it follows from the maximum principle that $w\le 0$ on $[0,T]\times [0,1]$ and we conclude.
\vs

In the argument above we ignored any boundary conditions for $\pt_p^2V$ without showing that they are irrelevant. There are two ways to  complete the proof. The first is to show, as we did earlier in the paper,  that there is no need for boundary conditions. The second is to show that  $\widetilde w=\zeta w\le 0$,  where $\zeta$ is a nonnegative function that vanishes near  $0$ and $1$. This implies that $w\le 0$  in the support of $\zeta$ and we may  conclude by letting the support of $\zeta$ become $(0,1)$. We leave the details to the interested reader.
%The second is to change variables to transform the problem to one defined on $[0,T]\times \R$ and to show that the transformed $w$ is nonnegative. We leave the details to the interested reader.

\end{proof}

\subsection{The geometry of the free boundaries}

We establish  some basic properties of the  
 exploration and stopping regions $ \cE^\mu$ and $\cS^\mu$, which are  defined in {\eqref{vlakas01} and \eqref{vlakas02}} respectively, and their $t-$slices are given by {\eqref{eq:two-boundary}} below. 
%\begin{equation}
% \begin{aligned}
% \cE^\mu&:=\{(t,p)\in[0,T)\times(0,1):V^\mu(t,p)<G(p)\},\\
%\cS^\mu&:=
%  \{(t,p)\in[0,T]\times[0,1]:V^\mu(t,p)=G(p)\},
% \end{aligned}
% \label{eq:exploration-stopping1}
%\end{equation}
%and, for $t<T$, we denote the spatial sections of $\cE^\mu$ and $\cS^\mu$ by
%\[
% \cE_t^\mu:=\{p\in(0,1):V^\mu(t,p)<G(p)\} \ \text{and} \ 
% \cS_t^\mu:=\{p\in[0,1]:V^\mu(t,p)=G(p)\}.
%\]

We  begin by showing  that the form of the obstacle $G$ yields that the exploration region contains a neighborhood of  the line $[0,T)\times \{m\}$. 
% when $G$ is piecewise affine and the whole region $\{F^\mu<0\}$ when $G$ is piecewise strictly concave. 

%the open set $\{F^\mu<0\}$. Since this set is empty when $G$ is piecewise linear, we  state a separate result for this case asserting that the exploration zone contains a neighborhood  of the line $[0,T]\times \{m\}$.
%% when $G$ is piece wise linear or the open set $\{F^\mu<0\}$ otherwise. Note that the latter set is empty when $G$ is piecewise linear. 
%\vs
%
%We begin with the latter case.

\begin{lemma}%[Negative forcing produces exploration: smooth case]
\label{lem:smooth-exploration}
Assume \eqref{takis20} and fix $\mu\in \cM$. Then, 
% If $G$ satisfies \eqref{ass:obstacle1}(i), 
\begin{equation}
 V^\mu(\cdot,m)<G(m) \ \ \text{in} \ \ [0,T).
 \label{eq:kink-exploration}
\end{equation}
%and, if $G$ satisfies \eqref{ass:obstacle1}(ii), then 
%\begin{equation}\label{eq:F-negative-exploration}
%\{(t,p)\in[0,T)\times(0,1): F^\mu(t,p)<0\}
% \subset\cE^\mu \ \ \text{and} \ \ \xi_-^\mu(t),\ \xi_+^\mu(t) \in\cE_t^\mu. 
% \end{equation}
 \end{lemma}
%and 
%%\begin{equation}
%% \xi_-^\mu(t),\ \xi_+^\mu(t)\in\cE_t^\mu,
%% \label{eq:forcing-zeros-exploration}
%%\end{equation}
%and there is $\delta_{\mu,t}>0$ such that
%\begin{equation}
% \{p\in[0,1]:F^\mu(t,p)<\delta_{\mu,t}\}
% \subset\cE_t^\mu.
% \label{eq:small-forcing-exploration}
%\end{equation}
%\end{lemma}

%\end{document}

\begin{proof}
To prove \eqref{eq:kink-exploration}, we fix $t\in [0,T)$, and,   for $\eps>0$ small, we choose $p^\ep_{ \pm} \in (0,1)$ such that 
\begin{equation}\label{elias1}
G(p^\ep_{\pm})=G(m)-(\overline C_f +1)\ep.
\nonumber
\end{equation}
Let $\overline u$ be the solution to 
\begin{equation}\label{elias2}
\begin{cases}
&-\pt_s \overline u-a^\mu (s,p) \pt_p^2 \overline u - \overline C_f=0 \ \ \text{in} \ \ [t,t+\ep) \times (p^\ep_{-}, p^\ep_{+}),\\[1.2mm]
&\overline u (\cdot,p^\ep_{\pm})=G(m)-(\overline C_f +1)\ep \ \ \text{in} \ \ [t,t+\ep],\\[1.2mm]
& \overline u(t+\ep,\cdot)=G \ \ \text{on} \ \ [p^\ep_{-}, p^\ep_{+}].
\end{cases}
\end{equation}
Since $V^\mu$ is a sub-solution to \eqref{elias2}, it follows that 
\beq\label{elias3}
V^\mu \leq \overline u \ \ \text{in} \ \  [t,t+\ep) \times [p^\ep_{-}, p^\ep_{+}].
\eeq
Moreover, 
\beq\label{elias4}
\overline u(s,p)\leq \big(G(m) + \overline C_f (t+\ep -s)\big)u_f(s,p) + \big(G(p^\ep_{\pm}) + \overline C_f (t+\ep -s)\big) u_b(s,p) \ \ \text{in} \ \ [t,t+\ep) \times [p^\ep_{-}, p^\ep_{+}],
\eeq
where  $u_b$ and $u_f$ solve the backward boundary value problems
\beq\label{takis40}
 -\pt_s u-a^\mu(s,p)\pt^2_p u=0 \ \  \text{in} \ \ [t,t+\ep)\times (p^\ep_{-}, p^\ep_{+}), 
\eeq
with 
\[u_b(\cdot, p^\ep_{\pm})=1 \ \ \text{ on} \ \  [t,t+\ep) \ \ \text{and} \ \  u_b(t+\eps,\cdot)=0 \ \ \text{on} \ \ [p^\ep_{-}, p^\ep_{+}],\]
and
\[u_f(\cdot,p^\ep_{\pm})=0  \ \ \text{ on} \ \   [t,t+\ep)  \ \ \text{and} \ \ u_f(t+\eps,\cdot)=1 \ \ \text{on} \ \ [p^\ep_{-}, p^\ep_{+}].\]

Note that the boundary conditions and the linearity of \eqref{takis40} imply that 
\beq\label{takis41}
u_b+u_f=1 \ \  \text{in} \ \ [t,t+\ep]\times [p^\ep_{-}, p^\ep_{+}].
\eeq

Using \eqref{elias3}, \eqref{elias4} and \eqref{takis41} we deduce that 
\beq\label{elias5}
V(t,m)-G(m)\le \ep\big( (\overline C_f+1) u_f(t,m) -1\big).
\nonumber
\eeq

The uniform ellipticity of \eqref{takis40} and elementary estimates about exit times imply that, as $\ep\to 0$, $u_f(t,m)\to 0$, and, hence, \eqref{eq:kink-exploration} holds.
\vs

\end{proof}

The next  result establishes  that the stopping region contains a, uniform in $\mu$, compact neighborhood of $[0,T]\times \{0,1\}$. The proof is based on the construction of uniform in $\mu$ explicit barriers near the end points of $[0,1]$. 

\begin{lemma}%[Uniform stopping neighborhoods of the endpoints]
\label{lem:endpoint-neighborhoods}
Assume \eqref{takis20}. Then, there exists $\kappa>0$, which is  independent of $(t,\mu)\in [0,T]\times \cM$, such that
\begin{equation}
 V^\mu=G
 %\quad\text{in }p\in[0,\kappa]\cup[1-\kappa,1] 
 \quad\text{in } \quad [0,T]\times \big([0,\kappa]\cup[1-\kappa,1]\big)\times \cM.
 \label{eq:endpoint-stopping-neighborhoods}
\end{equation}

\end{lemma}

\begin{proof}
%In view of \eqref{eq:kink-exploration} we divide the domain in two parts, namely $[0,T]\times [0,m]$ and $[0,T]\times [m,1]$.
%Since the arguments are similar,  here we present the details for the latter. 
%
%\vs

Let 
\beq\label{takis45}
M:=\frac{{\underline c_f}}{48\overline c_b^2},
%K_{G}:=\max\Big\{\underset{p\in [0,m]} \sup |G''(p), |\underset{p\in [m,1]} \sup |G''(p)|\Big\} \ \ \text{ and} \ \  M:=\frac{\underline f}{48\overline c_b^2},
\eeq
%In the smooth case set $K_G:=\|G''\|_\infty$; in the piecewise-affine case
%set $K_G:=0$, since $G$ is affine near the endpoints.  Put
%\[
% M:=\frac{\underline f}{48\overline\gamma^2}.
%\]
and choose $\kappa>0$ small enough so that 
%$4\kappa<1$ and, in the
%piecewise-affine case, $4\kappa<\min\{m,1-m\}$, and such that
\begin{equation}
 \frac{3M}{\kappa}>2\ol a \ \ \text{and} \ \ 2\kappa <m < 1-2\kappa,
 % \ \ \text{and} \ \ 
 %2\overline\gamma^2\kappa^2K_G\le\frac{\underline f}{4},
 \label{eq:kappa-choice}
\end{equation}
Let
\begin{equation}
 d_\kappa:=\frac{G(1-2\kappa)-G(2\kappa)}{1-4\kappa} \ \ \text{and} \ \ 
 L_\kappa(p):=G(2\kappa)-M+d_\kappa(p-2\kappa),
 \label{eq:asymmetric-middle-chord}
\nonumber
\end{equation}
and consider  the continuous function
\begin{equation}
 \underline v(p):=
 \begin{cases}
 G(p)-M\bigl(\dfrac{p}{\kappa}-1\bigr)_+^3
       &\text{in} \ \ [0,2\kappa],\\
 L_\kappa(p)&\text{in} \ \ [2\kappa,1-2\kappa],\\
 G(p)-M\bigl(\dfrac{1-p}{\kappa}-1\bigr)_+^3
       &\text{in} \ \ [1-2\kappa,1].
 \end{cases}
 \label{eq:endpoint-barrier}
\nonumber
\end{equation}
Using that  the chord joining
$(2\kappa,G(2\kappa))$ and
$(1-2\kappa,G(1-2\kappa))$ lies below $G$, we deduce that 
$L_\kappa\le G-M<G \ \text{ on} \  [2\kappa,1-2\kappa],$
while on the two remaining
pieces, it is immediate that 
$\underline v\le G.$
%Thus $\underline v\le G$ without
%any symmetry assumption.  
\vs

Direct calculations yield 
%$(\kappa,2\kappa)$,
\[
 \underline v''(p)=G''(p)-\frac{6M}{\kappa^2}\left(\frac p\kappa-1\right) \ \ \text{in} \ \ (\kappa,2\kappa),
\]
and 
\[
 \underline v''(p)=G''(p){-}\frac{6M}{\kappa^2}\left(\dfrac {1-p}{\kappa}-1\right) \ \ \text{in} \ \ (1-2\kappa,1-\kappa).
\]
%with the convention $G''=0$ in the piecewise-affine case.  

Since
{$G''=0$ on these intervals and $a^\mu(t,p)\le2\overline c_b^2\kappa^2$}, \eqref{takis45} and 
\eqref{eq:kappa-choice} give 
\begin{equation}
 {-a^\mu\underline v''-f^\mu
 \le
 12\overline c_b^2M-\underline c_f
 =-\frac34\underline c_f<0} \ \ \text{in} \ \ (\kappa,2\kappa)\cup (1-2\kappa,1-\kappa). 
 \label{eq:left-barrier-PDE}
\nonumber
\end{equation}
%The same calculation holds on $(1-2\kappa,1-\kappa)$.  
In the middle
interval $(2\kappa,1-2\kappa)$,  it is immediate that 
\[
 -a^\mu\underline v''-f^\mu=-f^\mu\le-{\underline c_f},
\]
while, on $[0,\kappa]\cup[1-\kappa,1]$, we have 
\[
 -a^\mu\underline v''-f^\mu
 {=-f^\mu\le-\underline c_f<0.}
\]
At the two junctions, that is, at $p=2\kappa$ and $p=1-2\kappa$,  using \eqref{eq:kappa-choice} and the fact that 
$|d_\kappa|,|G'(2\kappa)|,|G'(1-2\kappa)|\le L_G$, we find that 
\begin{align}
 \pt_p^+\underline v(2\kappa)
 -\pt_p^-\underline v(2\kappa)
 &=d_\kappa-G'(2\kappa)+\frac{3M}{\kappa}>0 \ \ \text{and}  \notag\\
 \pt_p^+\underline v(1-2\kappa)
 -\pt_p^-\underline v(1-2\kappa)
 &=G'(1-2\kappa)+\frac{3M}{\kappa}-d_\kappa>0.
 \label{eq:barrier-upward-jumps}
\nonumber
\end{align}
%where \eqref{eq:kappa-choice} and
%$|d_\kappa|,|G'(2\kappa)|,|G'(1-2\kappa)|\le L_G$ were used.  
Hence, there does not exist a 
$C^2-$test function which can touch $\underline v$ from above at either point.
\vs
It follows that  $\underline v$ is a viscosity sub-solution to  {\eqref{takis100}}, and,   since
$\underline v\le G$ at $t=T$, the comparison gives
\[\underline v\le V^\mu\le G.\] 

Because $\underline v=G$ on the two
endpoint intervals, \eqref{eq:endpoint-stopping-neighborhoods} follows.

\end{proof}

Next we recall the sets $A^\mu_\pm(t)$ given, for each $t\in [0,T)$, by \eqref{eq:boundary-definitions-intro}.

%Since, in view of 
%Lemma~\ref{lem:smooth-exploration} and 
%Lemma~\ref{lem:endpoint-neighborhoods}, for each $t\in [0,T)$, the set  $\cE_t^\mu$ is
%nonempty and bounded away from $0$ and $1$ and contains a neighborhood of $m$, it is immediate that
%\beq\label{eq:boundary-order}
%\kapp < A_-^\mu(t)<m< A_+^\mu(t)<1-\kappa. 
%\eeq
%endpoints by
%\begin{equation}
% A_-^\mu(t):=\inf\cE_t^\mu,\qquad
% A_+^\mu(t):=\sup\cE_t^\mu.
% \label{eq:boundary-definitions}
%\end{equation}
%Throughout the sequel, write
%\begin{equation}
% U^\mu(t,p):=G(p)-V^\mu(t,p)\ge0.
% \label{eq:U-definition}
%\end{equation}

The next proposition characterizes the $t-$slices  $\cE_t^\mu$ and $ \cS_t^\mu$ of the continuation and stopping regions in terms of $A_-^\mu(t)$ and  $A_+^\mu(t)$.

\begin{proposition}%[Two-boundary geometry]
\label{prop:two-boundary}
Assume \eqref{takis20} and fix $\mu\in \cM$. Then, there exists $\kappa\in (0,1)$ such that, for every $t \in [0,T)$,
\begin{equation}
 0<\kappa<A_-^\mu(t)<m<A_+^\mu(t)<1-\kappa<1. 
 \label{eq:boundary-order}
\nonumber
\end{equation}
Moreover, for each $t\in [0,T)$,
\begin{equation}
 \cE_t^\mu=(A_-^\mu(t),A_+^\mu(t)) \ \ \text{and} \ \ 
 \cS_t^\mu=[0,A_-^\mu(t)]\cup[A_+^\mu(t),1].
 \label{eq:two-boundary}
\end{equation}
%If $G$ satisfies \eqref{ass:obstacle1}(ii), then, for every $t \in [0,T)$,  %if there exist $\xi_\pm^\mu(t)$ as in \eqref{ass:obstacle1},
%\begin{equation}
% A_-^\mu(t)<\xi_-^\mu(t)<m<\xi_+^\mu(t)<A_+^\mu(t).
% \label{eq:strict-boundary-zero-separation}
%\end{equation}
%Finally, 
%\begin{equation}
% \cE_t^\mu=(A_-^\mu(t),A_+^\mu(t)) \ \ \text{and} \ \ 
% \cS_t^\mu=[0,A_-^\mu(t)]\cup[A_+^\mu(t),1].
% \label{eq:two-boundary}
%\end{equation}
%If there exist In both the piecewise-affine and smooth cases,
%$A_-^\mu(t)<m<A_+^\mu(t)$.
%In the smooth case, more precisely,
%\begin{equation}
% A_-^\mu(t)<\xi_-^\mu(t)<m<\xi_+^\mu(t)<A_+^\mu(t).
% \label{eq:strict-boundary-zero-separation}
%\end{equation}
\end{proposition}

\begin{proof} For each  $s\in [0,T)$, we use  Lemma~\ref{lem:endpoint-neighborhoods} and define the nearest contact points on the two sides of $m$ by
\beq\label{takis52}
 b_-^\mu(s)
 :=\max\{p\le m: V^\mu(s,p)=G(p)\} \ \ \text{and} \ \  b_+^\mu(s)
 :=\min\{p\ge m: V^\mu(s,p)=G(p)\},
\nonumber
\eeq
%
%\]
%and
%\[
% b_+^\mu(s)
% :=\min\{p\ge\m: V^\mu(s,p)=G(p)\}.
%\]
%These pointsexist by Lemma~\ref{lem:endpoint-neighborhoods}.
%  Their
It follows from the definition  above and Lemma~\ref{lem:smooth-exploration} that 
\[
 V^\mu(s,p)<G(p)  \ \ \text{if } \ \ p\in (b_-^\mu(s),b_+^\mu(s)).
\]

We prove that there is no continuation region to the right of
$b_+^\mu$ and left of $b_-^\mu$, which, in turn,  implies that, for all $s\in [0,T)$,  %the claim together with the fact that 
\[A_-^\mu(s)=b_-^\mu(s) \ \ \text{and} \ \ A_+^\mu(s)=b_+^\mu(s).\]
Since the arguments are similar,  we only present the details for the former. 
\vs

Fix $s\in[0,T)$ %write $b:=b_+^\mu(s)$, and set
and set $\wh V^\mu(s,p):=V^\mu(s,p) -G(p)$. It follows that   $\wh V^\mu\le 0$ and $\wh V^\mu(s,b_+^\mu(s))=\wh V^\mu(s,1)=0$. Moreover, in view of the linearity of $G$ on $[m,1]$,  Proposition~\ref{prop:concavity} yields that the map 
$p \to \wh V^\mu(s,p)$ is concave on $[m,1]$.
\vs
Hence, for all $p\in[b_+^\mu(s),1]$, 
\beq\label{takis53}
0{\ge} \wh V^\mu(s,p)\ge\frac{1-p}{1-b^\mu_+(s)}\wh V^\mu(s,b^\mu_+(s))+\frac{p-b^\mu_+(s)}{1-b^\mu_+(s)}V^\mu(s,1)=0. 
\nonumber
\eeq
Applying the same argument on $[0,m]$ and using that  $V^\mu(s,0)={\wh V^\mu}(s,b_-^\mu(s))=0$, we deduce that 
\beq\label{takis54}
V^\mu(s,p)=G(p)\quad\text{for }\quad p\in
[0,b_-^\mu(s)]\cup[b_+^\mu(s),1].
\nonumber
\eeq
This proves the required identities $A_\pm^\mu(s)=b_\pm^\mu(s)$.

\vs

The rest of the claim follows easily.

\end{proof}

\subsection{The stability of the free boundaries}

We state and prove  the stability of the free boundaries  with respect to the weak convergence of the collective input. 
The argument is complicated and the analysis {relies, among others, on two results} obtained in \cite{BDM2006} and  \cite{Blanchet2006exercise},
% Blanchet, Dolbeault and Monneau {\cite{BDM2006}} (see also Blanchet {\cite{Blanchet2006exercise}}), 
which are stated below. We note that these references 
%{\cite{BDM2006}} (and {\cite{Blanchet2006exercise}}) 
work with the forward in time problem for solutions to the obstacle problem with obstacle $0$, that is, they formulate their problem for $G-V^\mu$. 
\vs

Since our setting, motivated by the application,  was for $V^\mu-G$, below we rewrite the results of {\cite{BDM2006}}  in a form consistent with  the backward in time evolution. 
\vs

We begin with two important facts about the interfaces which we take from {\cite{BDM2006}} (see also {\cite{Blanchet2006exercise}}). To use them, we need to  reformulate the problem as we explain below. % of the variational inequality. 

{From} the discussion in the previous subsection, we know that the free boundaries  $\big(A_\pm^\mu (t)\big)_{t\in [0,T)}$ 
are in the domain  $ [0,T)\times \Big([0,1]\setminus\{m\}\Big)$. % \{(s,p)\in [0,T)\times [0,1]: F^\mu(s,p)>0\}.$
\vs
It is  a classical fact in the free boundary literature, that,  around the curves $\big(A_\pm^\mu(t)\big)_{t\in [0,T]}$
% and $\big(B^mu(t)\big)_{t\in[0,T)}$, 
we may  rewrite   {\eqref{takis100}} with unknown $\widehat V^\mu =V^\mu-G$ as 
\beq\label{v3}
\pt_t \widehat V^\mu +a^\mu \pt^2_p \widehat V^\mu =-f^\mu \mathbf{1}_{\{\widehat V^\mu<0\}} \ \ \text{and} \ \ \widehat V^\mu \le 0.
% \ \ \text{in a subset of $\{(s,p)\in [0,T)\times [0,1]: F^\mu(s,p)>0\}$.} %} \ \ Q_R(s_0,q_0),
\eeq
Indeed, 
%comparison with the supersolution $G+\varepsilon$ gives
%$0\le\beta_\varepsilon(V^{\mu,\varepsilon}-G)\le\beta(1)$.
standard  parabolic estimates and passage to the penalization limit
yield
$\widehat V^\mu\in W_{\mathrm{loc}}^{2,1;q}
((0,T)\times((0,1)\setminus\{m\}))$ for every $1<q<\infty$.
Since $\partial_t\widehat V^\mu$ and
$\partial_p^2\widehat V^\mu$ vanish almost everywhere on
$\{\widehat V^\mu=0\}$, the equation in the continuation region
gives \eqref{v3} almost everywhere.

\vs
In what follows, to state the results we use the parabolic cylinders around a point $(s_0,q_0)\in [0,T)\times (0,1)$ given by 
\[Q_R(s_0,q_0)=\{(s,q)\in [0,T)\times [0,1] : |s-s_0|<R^2 \ \text{and} \ |q-q_0|<R\},\]
and 
\[Q_R^+(s_0,q_0)=\{(s,q)\in [0,T)\times [0,1] : 0\le s-s_0<R^2 \ \text{and} \ |q-q_0|<R\}.\]

We also recall that our assumptions imply that there exist some $\delta>$ such that in a neighborhood of the free boundaries  we have 
%In the rest of this section and section~4 we will several times the following two results of \cite{BlDM}/ % and \cite{Bl}. 
%The key assumptions are that, for some $\delta>0$, 
\beq\label{ellie1}
a^\mu\ge \delta \ \ \text{and} \ \ f^\mu\ge \delta.
\nonumber
\eeq
%We remark  that \eqref{ellie1} holds in our setting near the interfaces.
%\vs

The first result is Theorem~2.1 in \cite{BDM2006}, as corrected by Theorem~1 in \cite{BDM2010}. We note that there is no $\mu$ in this reference. Nevertheless, the bound depends only on the regularity and bounds on the coefficients, which, in view of \eqref{takis20}, are uniform in $\mu$. 

\begin{proposition} \label{takis60.00} Given $(s_0,q_0) \in [0,T)\times \Big([0,1]\setminus\{m\}\Big)$ and for all $Q_R(s_0,q_0) \subset  [0,T)\times \Big([0,1]\setminus\{m\}\Big)$,  the $\wh V^\mu$'s are bounded in $W^{2,1;\infty}(Q_{R}(s_0,q_0))$ uniformly in $R\in (0,R_0)$ and $\mu\in \cM$.
\end{proposition}

The second result we need is the following non-degeneracy lemma ({Lemma~2.5 in \cite{BDM2006}}), which can be interpreted as the fact that the free boundaries are not ``blurred'', that is, they   cannot appear or disappear suddenly. We also note that, similarly to  the bounds in  {Proposition~\ref{takis60.00}}, the constants in the estimate below are independent of $\mu$. 
\vs
% For the statement we need to consider the forward cubes
%\[Q_R^+(s_0,q_0)=\{(s,q)\in [0,T]\times [0,1] : 0\le s-s_0<R^2 \ \text{and} \ |q-q_0|<R\}.\]
%

\begin{proposition}\label{takis51.00} Let   $\wh V^\mu$ be a solution to  \ \eqref{v3} in $Q^+_r(P)$ for $P=(s_0,p_0)\in \overline \{\wh V^\mu<0\}$ for sufficiently small $r>0$. Then,  there exist $\overline C$ and $\overline r$, which are independent of  $\mu\in \cM$ and $P$,  such that,  if $Q^+_r(P)\cap\{\wh V^\mu=0\}$ is nonempty, then, for all $r\in (0, \overline r)$,
\beq\label{takis65}
\underset{(s,q)\in Q^+_r(P)}\inf \wh V^\mu(s,q)\le -\overline C r^2.
\nonumber
\eeq
\end{proposition} 

%We also recall that \cite{Bl} yields a  set $\mathcal R\subset(0,T)$ of full measure on which the interfaces $A_\pm^\mu$ are H\"older continuous with exponent $1/2$. 
%\vs

Following \cite{Blanchet2006exercise}, we say that  the free boundary consists of ``regular'' and ``singular'' points. The set $\mathcal R$ of regular points has full measure and is locally a $C^{1/2}-$ graph. Moreover, around points of $\mathcal R$, the free boundary is described locally by a $C^{1/2}-$subgraph or by a $C^{1/2}-$supergraph.  Finally, if $(t,p)$ is a regular point, then 
\[\underset{(s,q)\to (t,p), (s,q)\in \{\wh V^\mu<0\}} \lim \pt_t \wh V^\mu(s,q)=0.\]

We {are} now ready for the stability result for the free boundaries with respect to the weak convergence of measures in $\cM$.

\begin{theorem}%[Stability of the stopping boundaries]
\label{thm:boundary-ae}
Assume \eqref{takis20} and let {$\mu_n,\mu\in\cM$} be such that, as $n\to \infty$, $\mu_n\weakto\mu$.  Then,  
there exists a set $\mathcal R\subset(0,T)$ of full
Lebesgue measure such that, for 
every compact $K\subset\mathcal R$,
\begin{equation}
 \lim_{n\to\infty}\sup_{t\in K}
 |A_\pm^{\mu_n}(t)-A_\pm^\mu(t)|=0.%+|B^{\mu_n}(t)-B^\mu(t)|\right)
 \label{eq:boundary-pointwise}
\nonumber
\end{equation}
In particular, for almost every $t\in(0,T)$, $A_\pm^{\mu_n}(t)\to A_\pm^\mu(t)$.
\end{theorem}

\begin{proof}
Since the arguments are similar, we only present the details for the convergence of the {$A_+^{\mu_n}$}'s to $A_+^\mu$.
\vs
Set $\wh V^n:=\wh V^{\mu_n}$, $\wh V:=\wh V^\mu$, {$A_\pm^n:=A_\pm^{\mu_n}$ and $A_\pm:=A_\pm^\mu$}, and recall that, in view of
Proposition~\ref{prop:value-stability}, we know that, as $n\to\infty$, {$\wh V^n$}$\to\wh V$ uniformly in $[0,T]\times[0,1]$.
\vs

Fix $t\in(0,T)$ and let $(t_n)_{n\in\N}$ be any sequence such that, as $n\to\infty$, $t_n\to t$. Since the $A_+^n(t_n)$'s are bounded,
along a subsequence, which, for notational simplicity, we denote as the sequence, and for some $A_+^\ast(t)\in[m,1-\kappa]$ (by Lemma~\ref{lem:endpoint-neighborhoods}), we have
\beq\label{chris1}
A_+^n(t_n)\to A_+^\ast(t).
\nonumber
\eeq

Recall that, if $p\in(A_-(t),A_+(t))$,  then
$\wh V(t,p)<0$. It then follows from the continuity of $\wh V$ and the uniform convergence
$\wh V^n\to\wh V$, that $\wh V^n(s,p)<0$ for all large $n$ and all $s$
sufficiently close to $t$. Proposition~\ref{prop:two-boundary} then
gives $p\in(A_-^n(s),A_+^n(s))$. Letting
$p\to A_+(t)$ yields
\begin{equation}
 A_+^\ast(t)\ge A_+(t).
 \label{eq:easy-boundary-limits}
\end{equation}

Next, Theorem~1.1 in \cite{BDM2006} 
implies that there exists  a set $\mathcal R\subset(0,T)$ of full measure such that
\[
 \partial_t\wh V(s,q)\longrightarrow0
 \quad\text{as }(s,q)\to(t,p)
\]
at every free boundary point $(t,p)$ with $t\in\mathcal R$
and $p\in(0,1)\setminus\{m\}$.

\vs
Fix $t\in\mathcal R$.
Since $A_+^\ast(t)>m$ by \eqref{eq:easy-boundary-limits},
Proposition~\ref{takis51.00} gives $\overline C$ and $\overline r$ such that, for every $\rho\in(0,\overline r)$ and $n$ sufficiently large,
\beq\label{vlakas2.1}
\inf_{\substack{0<s-t_n<\rho^2\\
 |q-A_+^n(t_n)|<\rho}}\wh V^n(s,q)\le-\overline C\rho^2.
\nonumber
\eeq
For each fixed $\rho$, these cylinders lie in
$Q_{2\rho}(t,A_+^\ast(t))$ for all large $n$.
Letting $n\to\infty$ and using uniform convergence gives, for every
sufficiently small $\rho>0$,
\beq\label{vlakas3}
 \inf_{Q_{2\rho}(t,A_+^\ast(t))}\wh V\le-\overline C\rho^2.
\eeq

This, in turn, implies that
\beq\label{vlakas4}
A_+^\ast(t)=A_+(t).
\nonumber
\eeq
Indeed, if $A_+^\ast(t)>A_+(t)$, then $\wh V(t,\cdot)=0$ on
a neighborhood of $A_+^\ast(t)$. The defining property of $\mathcal R$
and integration in time give
\[
 \sup_{Q_r(t,A_+^\ast(t))}|\wh V|
 \le r^2\operatorname*{ess\,sup}_{Q_r(t,A_+^\ast(t))}
 |\partial_t\wh V|=o(r^2).
\]
Taking $r=2\rho$ contradicts \eqref{vlakas3}.
\vs
The same argument for $A_-$ and the arbitrariness of the subsequence give
$A_\pm^{\mu_n}(t_n)\to A_\pm^\mu(t)$ whenever $t_n\to t\in\mathcal R$.
Taking $\mu_n\equiv\mu$ also gives continuity at such $t$. Then, compactness
yields uniform convergence on each compact $K\subset\mathcal R$.

\end{proof}

\

\subsection{The temporal continuity of the free boundaries} % at contact points for beliefs at $t=0$}

Theorem~\ref{thm:boundary-ae} gives the convergence of the boundaries on 
a full-measure set of deterministic times. However, this alone does not imply
regularity at the random first-contact time, a  question we address next.
\vs
%The continuity of the interfaces in one space dimension has been analyzed extensively; see, for example, \cite{BlDM} and 
%\cite{BL}. Loosely speaking, outside a set of times of measure zero, the interfaces are locally H\"older continuous. 
The analysis in \cite{BDM2006} and \cite{Blanchet2006exercise}   as well as subsequent references is based on a careful blow up in space-time argument around each interface point and the complete characterization of the limit problems. 
\vs

Here,  we  investigate the local continuity of the free boundaries at the contact points of the belief process that starts at $P_0$ at $t=0$. For this, we   use the blow-up argument and the results of \cite{BDM2006} and 
{\cite{Blanchet2006exercise}} to show that, at any contact point $0<\tau(\omega)< T$ and a.s. in $\omega$, the interfaces are $1/2-$H\"older continuous.  Then, we investigate separately the continuity of the $A_\pm^\mu$'s as $t\to 0$ and $t\to T$. 
\vs

All these  are crucial {to} establish the existence of the mean field equilibria in the next section, which is based on showing the continuity in probability of the best-response pair  $(\hat\tau_{0,p}^\mu,\hat d_{0,p}^\mu)$   defined
in \eqref{eq:best-response-pair}.

\vs
We begin with a technical lemma which is a modification of a result of Krylov {\cite{Krylov2003sqrt}}. 

\begin{lemma}
\label{lem:posterior-square-root}
Fix $\mu\in{\mathcal M}$, and let $P^\mu$ be the  belief process \eqref{0}. For every $[0,T]$-valued $\mathbb F^\mu$-stopping time $\sigma$ and almost surely
on $\{0<\sigma\le T,\ P^\mu_\sigma \in(0,1)\}$, there exist $C>0$ and a sequence 
$t_k\uparrow\sigma$ such that
\[
 |{P^\mu_{t_k}}-P^\mu_\sigma|\le C\sqrt{\sigma-t_k}.
\]
\end{lemma}

\begin{proof}
Since $P^\mu-P^\mu_0$ is a continuous local $\mathbb F^\mu$-martingale with quadratic
variation
\[
 \langle P^\mu\rangle_t
 =\int_0^t(\Delta b(s;\mu_{\le s}))^2(P^\mu_s)^2(1-{P^\mu_s})^2\,\dd s,
\]
the Dambis--Dubins--Schwarz theorem \cite[Section~3.4.B]{KaratzasShreve1991}  provides a Brownian motion
$\big(\widetilde W_q\big)_{q\ge0}$, on, if necessary,  enlarged probability space, such that, for all $t\in [0,T]$, 
\[
 P^\mu_t=P^\mu_0+\widetilde W_{\langle P^\mu\rangle_t}, 
% \qquad 0\le t\le T
\]
that is,  $P^\mu$ is a Brownian motion evaluated at its quadratic-variation
clock and, moreover,  $\langle P^\mu\rangle_T\le\overline c_b^2T/16$.
\vs
%Since $0\le [p(1-p)]^2\le1/16$ on $[0,1]$, we have 
%$\langle P\rangle_T\le\bar\gamma^2T/16$.
%\vs
It follows from equation~$(2.1)$ of Theorem~$2.1$ in \cite {Krylov2003sqrt} that, almost surely,
for every {$q\in(0,\overline c_b^2T/16]$}, there exist 
$r_k=2^{-j_k}\downarrow0$ {with $r_k<q$} such that, for some constant $c>0$,
\[
 \operatorname{osc}_{[q-r_k,q]}\widetilde W_.\le c\sqrt{r_k}.
\]
Let $q=\langle P^\mu\rangle_\sigma$. Since
$P^\mu_\sigma\in(0,1)$, the continuity of $P^\mu$ and the bound 
$\Delta b^\mu(s)\ge\underline c_b$ imply that the map $t\to \langle P^\mu\rangle_t$ is
strictly increasing on a left neighborhood of $\sigma$. Hence, for
all sufficiently large $k$, there exist $t_k\uparrow\sigma$ such that 
$\langle P^\mu\rangle_{t_k}=q-r_k$. % and $t_k\uparrow\sigma$.  
It then follows that  %Consequently,
\[
 |P_{t_k}-P_\sigma|
 \le c\sqrt{r_k} \ \ \text{with} \ \ 
 r_k=\langle P\rangle_\sigma-\langle P\rangle_{t_k}
 \le\frac{\overline c_b^2}{16}(\sigma-t_k),
\]
hence the claim.

\end{proof}

The next result is about the temporal continuity of the free boundaries at the contact points. The proof uses the blow-up analysis in Proposition~3.2, Theorem~3.9, Proposition~4.1 and Proposition~5.5 of \cite{BDM2006}, together with Proposition~2.2 and Theorem~1.1(i) of \cite{Blanchet2006exercise}, where we refer for some of the background. 

\begin{proposition}%[Local $C^{1/2}$ estimate at first contact]
\label{lem:local-contact-Holder}
Assume \eqref{takis20}, fix $\mu\in{\mathcal M}$, and let  $P^\mu:=P^{0,P_0;\mu}$ and 
$\tau:=\hat\tau_{0,P_0}^\mu$ with $\hat\tau_{0,P_0}^\mu$ given by \eqref{eq:best-response-pair}. Then, almost surely on
$\{0<\tau<T\}$, we have that  %the first contact occurs at one of the two free boundaries, 
\[
 P^\mu_\tau=A_-^\mu(\tau)\quad\text{or}\quad P^\mu_\tau=A_+^\mu(\tau).
\]
Moreover, there exist constants $C,r>0$, {which may depend on $\mu$ and the realized contact point}, such that, for all $s,s'\in(\tau-r,\tau+r)\cap(0,T)$, if $P^\mu_\tau=A_-^\mu(\tau)$, then
\[
 |A_-^\mu(s)-A_-^\mu(s')|\le C\sqrt{|s-s'|},
\]
and, if $P_\tau=A_+^\mu(\tau)$, then 
\[
 |A_+^\mu(s)-A_+^\mu(s')|\le C\sqrt{|s-s'|}.
\]
%whenever $s,s'\in(\tau-r,\tau+r)\cap(0,T)$.
\end{proposition}
%
%\begin{remark*}[Remarks on Lemma~\ref{lem:local-contact-H\"older}]
%The lemma proves two facts.  First, the realized first-contact point is
%a regular free-boundary point.  Following Blanchet \cite{{Blanchet2006exercise}}, this means that the
%coincidence set has positive lower density in backward parabolic
%cylinders; equivalently, in calendar time,
%\[
% \liminf_{r\downarrow0}
% \frac{\bigl|\{U^\mu=0\}\cap
% ((\tau,\tau+r^2)\times(P(\tau)-r,P(\tau)+r))\bigr|}
% {2r^3}>0.
%\]
%Second, \cite[Theorem~1.1(i)]{Blanchet2006exercise} turns this density
%property into a local $C^{1/2}$ free-boundary graph.  The resulting
%estimate for $A_-^\mu$ or $A_+^\mu$ is used in
%Proposition~\ref{prop:canonical-response-continuity} to prove convergence
%of the stopping times.
%\end{remark*}
%

\begin{proof}
Let $\Omega_0$ be a probability-one event carrying the
conclusions of Lemma~\ref{lem:posterior-square-root} with $\sigma=\tau$ and the continuity of the 
posterior process. In what follows, we work pathwise for each  fixed  $\ol \omega\in\Omega_0\cap\{0<\tau<T\}$.
% and work pathwise on .
%continuity, and invariance of $(0,1)$.
\vs
%Fix $\omega\in\Omega_0\cap\{0<\tau<T\}$ and work pathwise.
We set $p_0:=P^\mu_\tau\in \cS_\tau^\mu\cap(0,1)$.  Proposition~\ref{prop:two-boundary} then implies that 
$(\tau,p_0)$ is a space-time free-boundary
point with $p_0\ne m$. % by Proposition~\ref{prop:two-boundary}.}
\vs
We recall that,
as discussed earlier, near the free boundaries, $\wh V^\mu$ 
solves, in a small parabolic cubes $Q_R$ containing $(\tau,p_0)$,  the problem 
\beq\label{pl4}
 \pt_t {\wh V}^\mu+a^\mu \pt^2_p \wh V^\mu
 =
  -f^\mu{\1_{\{\wh V^\mu<0\}}} \ \ \text{a.e. on } \ \ Q_R,
 \eeq
 with $a^\mu$ and $f^\mu$ bounded from below away from zero.
\vs
%in view of the assumptions on the data, %there exists $k_0>0$ such that 
%\beq\label{pl1} k_0:=F^\mu(\tau,p_0)>0.\eeq
%As discussed earlier, near the free boundaries, $\wh V^\mu$ 
%solves, in a small parabolic cubes $Q_R$ containing $(\tau,p_0)$,  the problem 
%\beq\label{pl4}
% \pt_t {\wh V}^\mu+A_-^\mu \pt^2_p \wh V^\mu
% =
%  -F^\mu\1_{\{V^\mu<0\}} \ \ \text{a.e. on } \ \ Q_R,
% \eeq
% with $A_-^\mu$ and $F^\mu$ bounded from below away from zero.
%\vs

%By \eqref{eq:small-forcing-exploration} and
%\eqref{eq:coefficient-bounds}, $k_0>0$.  Moreover, in the
%piecewise-affine case, $p_0\ne m$ by
%Proposition~\ref{prop:two-boundary}.  Hence there are $R,\lambda>0$
%such that $Q_R:=(\tau-R^2,\tau+R^2)\times(p_0-R,p_0+R)
% \Subset(0,T)\times(0,1)$;
%$Q_R$ does not meet $p=m$ in the piecewise-affine case, and on $Q_R$
%the diffusion coefficient and the relevant right-hand side are
%H\"older continuous and bounded below by $\lambda$:
%\[
% A_-^\mu\ge\lambda,\qquad
% F^\mu\ge\lambda\ \text{in the smooth case},\qquad
% f^\mu\ge\lambda\ \text{in the piecewise-affine case}.
%\]
%Then Lemma~\ref{takis} yields the uniform bound 
% %Standard local regularity for parabolic obstacle problems gives
%\[
% U_t^\mu+A_-^\mu U_{pp}^\mu
% =
% \begin{cases}
%  F^\mu\1_{\{U^\mu>0\}},&\text{in the smooth case},\\
%  f^\mu\1_{\{U^\mu>0\}},&\text{in the piecewise-affine case},
% \end{cases}
% \qquad\text{a.e. on }Q_R.
%\]
%Then the estimate in \cite[Theorem~1]{BDM2010} then gives
Choose $R_0>0$ with
$Q_{2R_0}(\tau,p_0) \subset (0,T)\times\big((0,1)\setminus\{m\}\big)$.
Proposition~\ref{takis60.00} gives a constant $C_{R_0}>0$ such that 
\beq\label{takis72}
 \|\wh V^\mu\|_{W^{2,1;\infty}(Q_{R_0}(\tau,p_0))}\le C_{R_0}.
\eeq
Following  \cite{BDM2006} and \cite{Blanchet2006exercise}, we perform a blow-up 
argument near the point $(\tau,p_0)$.  We set  $a_0^\mu:=a^\mu(\tau,p_0)$
and $f^\mu_0:=f^\mu(\tau,p_0)>0$, and consider the family 
\[
 \wh V^{\mu,\eps}({s,p}):=\frac{1}{\ep^2}\wh V^\mu\!\left(
 \tau+\frac{\eps^2}{f^\mu_0}s,\,
 p_0+\eps\sqrt{\frac{a^\mu_0}{f^\mu_0}}\,p\right).
\]
%with $a_0:=A_-^\mu(\tau,p_0)$
In view of \eqref{takis72} and \cite[Proposition~3.2]{BDM2006},
every sequence $\eps\to 0$ admits a subsequence along which the 
$\wh V^{\mu,\eps}$'s converge, locally uniformly in $\R\times\R$, 
to a solution $\wh V^{\mu,0}:\R\times\R\to(-\infty,0]$ to 
\beq\label{thaleia100}
\pt_s\wh V^{\mu,0}+\pt_p^2\wh V^{\mu,0}
=-\1_{\{\wh V^{\mu,0}<0\}}
\quad\text{a.e. in }\R\times\R,
\nonumber
\eeq
and
\beq\label{thaleia101}
(0,0) \in \pt\{\wh V^{\mu,0}=0\}.
\nonumber
\eeq
These blow-up limits are classified by
Proposition~3.8, Theorem~3.9, and Proposition~4.1 in \cite{BDM2006},
after changing sign and reversing time.
\vs

It follows that all possible limits  have either one of the ``regular'' profiles 
\begin{eqnarray}\label{eq:regular_profiles}
     {\wh V^{\mu,0}}(s, p)=-\frac12(\max\{p,0\})^2 \ \ \text{or} \ \ 
 {\wh V^{\mu,0}}(s,p)=-\frac12(\max\{-p,0\})^2,
\end{eqnarray}
or are singular, in which case, for some $c\in [0,1/2]$, 
\[
 {\wh V^{\mu,0}}({0,p})=-cp^2. %\qquad 0\le c\le\frac12.
\]

%which, in view of 
%%Since $U^\mu(t,P(t))>0$ for $t<\tau$ and
%%$U^\mu(\tau,p_0)=0$, $(\tau,p_0)$ is a free-boundary point.
%%The preceding $W_{\rm loc}^{2,1;\infty}$ estimate and
%%\cite[Proposition~3.2]{BDM2006} imply that every sequence
%%$\eps_n\downarrow0$ has a locally uniformly convergent subsequence.
%Then  \cite[Proposition~4.1, Theorem~3.9, and
%Proposition~5.5]{BDM2006} yield that we  have either one of %yield that, every blow-up limit $u_0$ has  either one of
%the regular profiles
%\begin{eqnarray}\label{eq:regular_profiles}
%     v_0(x,s)=\frac12(\max\{x,0\})^2 \ \ \text{or} \ \ 
% v_0(x,s)=\frac12(\max\{-x,0\})^2,
%\end{eqnarray}
%or it is singular, in which case, for some $c\in [0,1/2]$, 
%\[
% v_0(x,0)=cx^2. %\qquad 0\le c\le\frac12.
%\]
%Since $p_0\in \cS^\mu_\tau$, we have either 
%If $p_0\in A_-^\mu(\tau)$ or $p_0\in A_+^\mu(\tau)$. 
Since $p_0\in\cS_\tau^\mu$, either $p_0\le A_-^\mu(\tau)$
or $p_0\ge A_+^\mu(\tau)$.
If $p_0\le A_-^\mu(\tau)$, for each fixed $p\le0$ and
$\ep$ small enough, we have \ 
$p_0+\eps\sqrt{\frac{a^\mu_0}{f^\mu_0}}\,p \le p_0\leq A_-^\mu(\tau)$
and, hence,
\[\wh V^{\mu,\ep}(0,p)=\frac{1}{\ep^2}\wh V^\mu(\tau, p_0+\eps\sqrt{\frac{a_0^\mu}{f^\mu_0}}\,p)=0,\]
which implies that 
\[\wh V^{\mu,0}(0,p)=0\quad\text{if }\quad p\le0.\]
If $p_0\ge A_+^\mu(\tau)$, for {each fixed} $p\ge0$ and $\ep$ small enough, we have \ 
$p_0+\eps\sqrt{\frac{a_0^\mu}{{f^\mu_0}}}\,p \ge p_0 \ge A_+^\mu(\tau)$
and, thus, 
\[\wh V^{\mu,\ep}(0,p)=\frac{1}{\ep^2}\wh V^\mu(\tau, p_0+\eps\sqrt{{\frac{a_0^\mu}{f^\mu_0}}}\,p)=0
\quad\text{and}\quad {\wh V^{\mu,0}(0,p)=0}\quad\text{if } \quad p\ge0.\]
%and
%\[u_0(0,x)\ge 0 \ \ \ \text{if}  \ \ x\ge 0.\]
%%It follows that in the first case $u_0(0,x)\le 0$ for $x\leq 0$, and in the second case $u(0,x)=0$ for $x\geq 0$. 
%\vs
As a result, in the singular case we must have $c=0$,
{so the classification leaves only the temporal solution}
\beq\label{takis80}
 {\wh V^{\mu,0}(s,p)}=\min\{-s,0\}=-s^+.
\nonumber
\eeq

We now rule out this limit. At a singular contact point, Proposition~5.5 in \cite{BDM2006} gives the same limit along every
blow-up sequence.
Indeed, applying Lemma~\ref{lem:posterior-square-root} with $\sigma=\tau$
and using the definition of $\tau$, we can find
 $t_k\uparrow\tau$ such that
\beq\label{pl2}
|P^\mu_{t_k}-p_0|\le C\sqrt{\tau-t_k} \ \ \text{and} \ \ \wh V^\mu(t_k,P^\mu_{t_k})<0.
\eeq
%while, in view of the definition of $\tau$,  $t_k<\tau$ implies that  $U^\mu(t_k,P(t_k))>0$.
\vs
Set 
\[\eps_k^2:={f_0}(\tau-t_k)\quad\text{and}\quad
p_k:=\frac{P^\mu(t_k)-p_0}{\sqrt{{a_0^\mu}(\tau-t_k)}}.\]
Then, in view of \eqref{pl2}, we know that the $p_k$'s are  bounded and, after passing, if necessary,  to a subsequence, we may assume
that $p_k\to p$, while we have $\wh V^{\mu, \eps_k}\to \wh V^{\mu,0}$ locally uniformly and 
\beq\label{pl3}
\wh V^{\mu,\eps_k}({-1,p_k})<0.\eeq
%Since $(x_k)$ is bounded, after passing to a subsequence we may assume
%that $x_k\to x$ and $u_{\eps_k}\to u_0$ locally uniformly.  
%\vs
%If the blow-up limit is singular, the preceding classification gives
%$u_0(x,s)=\max\{s,0\}$.  
%\vs
%Next, we recall that   \cite[Proposition~3.3]{BDM2006} that, as $k\to \infty$,  $\mathbf 1_{\{u_k>0\}} \to \mathbf{1}_{\{u>0\}}$ and $\mathbf 1_{\{u_k=0\}} \to \mathbf{1}_{\{u=0\}}$ almost everywhere.
%\vs
%It follows that, 
%$\max\{s,0\}=0$ near $(x,-1)$, 
%%\cite[Proposition~3.3]{BDM2006} implies that 
%the $u_{\eps_k}$'s  must have  a zero
%near $(x_k,-1)$ for all large $k$, while Lemma~\ref{takis51} (see also  
%\cite[Lemma~2.1]{Blanchet2006exercise}) gives fixed
%$0<\overline r$ and $\overline C>0$ such that, for all $r\le \overline r$,
%\[
%\overline C\rho^2
% \le
% \sup_{\substack{|y-x_k|<r\\-1<s<-1+r^2}}
% u_{\eps_k}(y,s)
%% \longrightarrow0,
%\]
%which is a contradiction, since, in view of the previous, $u_{\ep_k}\to 0$ in a neighborhood of $(x,-1)$.
%\vs 
%It follows that  $ p_0=A_-^\mu(\tau)$ or $ p_0=A_+^\mu(\tau)$.
%Thus, only the two limits in \eqref{eq:regular_profiles} remain, and, in view of %Then 
%\cite[Proposition~2.2]{Blanchet2006exercise},  
%$(\tau,p_0)$ is regular.
%\vs

The rescaled coefficients have uniform local H\"older and
positive lower bounds, so the nondegeneracy constants are uniform in
all sufficiently large $k$.
\vs
Choose $\delta\in(0,\overline r\wedge1/2)$, with $\overline r$ given in Proposition~\ref{takis51.00}
{for the rescaled equations},
% smaller than the radius allowed in the
%nondegeneracy estimate, 
and set
\[
 D:=\Big\{(s,y):
 |y-{p}|<\frac{\delta}{2},\
 {0<s+1<\frac{\delta^2}{4}}\Big\}.
\]
Since {$\wh V^{\mu,0}(s,y)=-s^+=0$} in $D$, Proposition~$3.3$ in 
{\cite{BDM2006}} and dominated convergence give that, as $k\to\infty$,
\[
 \bigl|D\cap\{\wh V^{\mu, \eps_k}<0\}\bigr|
 =
 \int_D\mathbf 1_{\{\wh V^{\mu, \eps_k}<0\}}(s,q) \dd q\dd s
 \longrightarrow0.
\]
{In view of the facts that $|D|>0$ and $\wh V^{\mu,\eps_k}\le0$, we then deduce that,
for all sufficiently large $k$,}
\[D\cap\{\wh V^{\mu, \eps_k}=0\}\ne\varnothing.\]
Moreover, since ${p_k\to p}$,  for all sufficiently large
$k$, every point $(s,y)\in D$ satisfies
$|y-{p_k}|<\delta$ and {$0<s+1<\delta^2$}, and, thus,
{$\wh V^{\mu,\eps_k}$ has a zero in the forward cylinder
$\{|y-{p_k}|<\delta,\ -1<s<-1+\delta^2\}$.}
\vs

It follows from \eqref{pl3} and Proposition~\ref{takis51.00}
(see \cite[Lemma~2.1]{Blanchet2006exercise}) that
{there exists a constant $\overline C>0$, which is  independent of all sufficiently
large $k$, such that, for this fixed $\delta$,}

\[
 \inf_{\substack{|y-{p_k}|<\delta\\-1<s<-1+\delta^2}}
 \wh V^{\mu,\eps_k}(s,y)\le-\overline C\delta^2.
\]

Because ${p_k\to p}$ and {$\delta<1/2$}, these cylinders lie in a fixed
compact subset of $\{s<0\}$, where $\wh V^{\mu,0}=0$. Hence,  the local uniform
convergence $\wh V^{\mu, \eps_k} \to \wh V^{\mu,0}$ gives
\[
 \sup_{\substack{|y-{p_k}|<\delta\\-1<s<-1+\delta^2}}
 |\wh V^{\mu,\eps_k}(s,y)|\longrightarrow0,
\]
which is a contradiction.

\vs 
{In conclusion, only the two regular limits in \eqref{eq:regular_profiles}
remain. But then, the strict inequality $p_0<A_-^\mu(\tau)$ or $p_0>A_+^\mu(\tau)$
would force an identically zero spatial trace, which  contradicts these
profiles. Hence $p_0=A_-^\mu(\tau)$ or
$p_0=A_+^\mu(\tau)$} Then, in view of
\cite[Proposition~2.2]{Blanchet2006exercise},  
$(\tau,p_0)$ is a regular point.
\vs

Next, we apply Theorem~$1.1(i)$  in \cite{Blanchet2006exercise} to \eqref{pl4}, the hypotheses of which are satisfied. Indeed,  our assumptions guarantee that the diffusion coefficients
and  $f^\mu$ in \eqref{pl4}
are H\"older continuous and bounded below away from $0$ on
the cylinder $Q_{R_0}(\tau,p_0)$ used in \eqref{takis72},
and \eqref{takis72} gives
$\wh V^\mu\in W_{\mathrm{loc}}^{2,1;\infty}$.
%Thus the hypotheses of Theorem~1.1(i) in
%\cite{Blanchet2006exercise} hold.
\vs

It follows that, since $(\tau,p_0)$ is a regular point,  there exist $\rho\in(0,R_0)$ and a function
\[
 g\in C^{1/2}((\tau-\rho^2,\tau+\rho^2)),\qquad g(\tau)=p_0,
\]
such that, in
$Q_\rho(\tau,p_0)\Subset(0,T)\times((0,1)\setminus\{m\})$,
the stopping set $\{\wh V^\mu=0\}$ is either
$\{p\le g(s)\}$ or $\{p\ge g(s)\}$. Finally,  we use  the continuity of $g$ to choose  $0<r<\rho^2$ so that
$|g(s)-p_0|<\rho$ whenever $|s-\tau|<r$.
\vs

Next, we identify $g$ as $A_-^\mu$ or $A_+^\mu$, depending on  which
boundary is hit.
\vs

Next, suppose that $p_0=A_-^\mu(\tau)$.
Then the  cylinder $Q_\rho(\tau,p_0)$ lies below $m$.
In view of  Proposition~\ref{prop:two-boundary}, at time $s=\tau$ and for
$|p-p_0|<\rho$, stopping is optimal exactly when
$p\le p_0=A_-^\mu(\tau)$. Hence, in this neighborhood,
\[
 \wh V^\mu(s,p)=0\quad\Longleftrightarrow\quad p\le g(s),
 \qquad |s-\tau|<r,\quad |p-p_0|<\rho.
\]
Since $g(s)$ is a stopping point smaller than  $m$, $g(s)\le A_-^\mu(s)$.
Then,  every $p\in(g(s),p_0+\rho)$ is a continuation point, so
$A_-^\mu(s)<p$. Letting $p\to g(s)$ gives
$A_-^\mu(s)\le g(s)$.
Therefore, for $s\in(\tau-r,\tau+r)$,
\[
 A_-^\mu(s)=g(s),
\]
and, thus,
\[
 |A_-^\mu(s)-A_-^\mu(s')|
 \le [g]_{C^{1/2}}|s-s'|^{1/2}.
\]
\vs

The case $p_0=A_+^\mu(\tau)$ follows similarly, yielding 
$A_+^\mu(s)=g(s)$ for $s\in(\tau-r,\tau+r)$ and the same
H\"older estimate.
\vs
Taking $C=1+[g]_{C^{1/2}((\tau-r,\tau+r))}$,
with $r$ chosen above, we conclude the proof.

\end{proof}

\begin{remark} Blanchet's theorem in \cite{Blanchet2006exercise} does not show that the entire free boundary is regular. The  singular profile {$\wh V^{\mu,0}(s,p)=-s^+$} may occur and the theorem gives local $C^{0,1/2}-$regularity only at points already known to be regular. The first contact point has an additional property which does not hold  at an arbitrary free boundary point. Indeed, {Lemma~\ref{lem:posterior-square-root}} provides continuation points approaching the free boundary point at a parabolic scale, which become, after rescaling,  bounded points $(-1, p_k)$ satisfying {$\wh V^{\mu,\eps_k}(-1, p_k)<0$}. This is, however, incompatible with {$-s^+$}, which vanishes for $s<0$. 
%by the phase convergence and the nondegeneracy estimate. 
Since the two boundary geometry excludes the remaining singular profiles, the first contact point is regular, and Blanchet's theorem then yields the local $C^{1/2}$ boundary graph.
\end{remark}

%\medskip
% \underline{Step 2: Local $C^{1/2}$ regularity.}
%The conditions already verified on $Q_R$ are precisely those required

Next, we turn to the continuity at $t=0$ which is necessary {to} examine what happens when $\tau=0$.

\begin{proposition}
%[Initial-time continuity and time-zero stability:
%piecewise-affine case]
\label{lem:initial-piecewise}
For every $\mu\in{\mathcal M}$,
\[
\underset{t\to 0} \lim A_\pm^\mu(t)=A_\pm^\mu(0).\]
% 
%Moreover, if 
%$\mu_n\rightharpoonup\mu$, then, with the notation of
%Proposition~4.1,
%\[
% \limsup_{n\to\infty}\tau_n\le0=\tau
% \qquad\text{almost surely on }\{\tau=0\}.
%\]
\end{proposition}

\begin{proof}
The proofs are similar for $A^\mu_-$ and $A^\mu_+$, we  only show the details  for the former.
\vs

Suppose that
$t_n\to 0$  and $A^\mu_-(t_n)\to a$. Since it is  immediate that $a\le A^\mu_-(0)$,  we assume  that $a<A^\mu(0).$  
\vs
Then, for every
$z\in(a,A_-^\mu(0))$, we know that  $\wh V^\mu(0,\cdot)=0$ in a neighborhood $I_z$ in $[0,1]$ of $z$, while $\wh V^\mu(t_n,z)<0$ for all large $n$.
\vs
Moreover, since, for some small $h>0$,
\[\partial_t \wh V^\mu + a^\mu \pt_p^2 \wh V^\mu=-f^\mu \ \ \text{in} \ \ [0,h)\times I_z,\]
it follows that 
\beq\label{eq:t_limit}
\partial^+_t  \wh V^\mu(0,p)=-f^\mu \ \ \text{in} \ \ I_z,
\eeq
where 
\[\partial^+_t  \wh V^\mu(0,p)=\underset{t \to 0+} \lim \dfrac{\wh V^\mu(t,p)}{t}.\]

%Choose $h>0$ such that $(-h,h)\times I_z$ does not meet $p=m$, and
%define
%\[
% \bar U(t,p):=
% \begin{cases}
%  0,&-h<t\le0,\\
%  U(t,p),&0<t<h.
% \end{cases}
%\]
%Extend $a$ and $f$ to $t<0$ by their values at $t=0$. Since
%$U(0,\cdot)=0$ on $I_z$, the gluing introduces no time-interface
%measure, and locally
%\[
% \bar U_t+\bar a\,\bar U_{pp}
% =\bar f\,\mathbf 1_{\{\bar U>0\}}
% \qquad\text{on }(-h,h)\times I_z.
%\]
%Extend \(U\) by zero to \(t<0\) on \(I_z\), and extend \(a\) and \(f\)
%by their values at \(t=0\). Since \(U(0,\cdot)=0\) on \(I_z\), the
%extension solves the obstacle equation locally across \(t=0\), so
%\((0,z)\) is an interior free-boundary point. Applying the blow-up
%argument of Step~1 in the proof of Lemma~3.15, the zero time-zero trace
%leaves only the temporal profile \(s^+\). Consequently,
%\begin{eqnarray}
%     \frac{U(t,z)}{t}\to f(0,z)\,\,
%\textrm{ as } \,\, t \to 0.
%\label{eq:t_limit}
%\end{eqnarray}
%
\vs
Now choose $a<p<q<A_-^\mu(0)$. Then, the concavity of the map $p\to \wh V^\mu({t_n},p)$ and the fact that 
$\wh V^\mu(t_n,A_-^\mu(t_n))=0$ give
\[
 \frac{\wh V^\mu(t_n,p)}{p-A_-^\mu(t_n)}
 \ge\frac{\wh V^\mu(t_n,q)}{q-A^\mu_-(t_n)}.
\]
Using \eqref{eq:t_limit}, {dividing by $t_n$}, and letting $n\to\infty$ yields
\[
 \frac{f^\mu(0,p)}{p-a}\le\frac{f^\mu(0,q)}{q-a},
\]
{a contradiction to  the concavity and strict positivity of
$f^\mu(0,\cdot)$, since
\[
 f^\mu(0,p)\ge
 \frac{q-p}{q-a}f^\mu(0,a)+\frac{p-a}{q-a}f^\mu(0,q).
 >\frac{p-a}{q-a}f^\mu(0,q), 
\]
Therefore, it must be that $\liminf_{t\downarrow0}A_-^\mu(t)\ge A_-^\mu(0)$.
\vs
The reverse upper bound follows from the continuity of $\wh V^\mu$ and
Proposition~\ref{prop:two-boundary}, completing the proof for $A_-^\mu$.}
 
 \end{proof}
 
To complete the results for the temporal  continuity properties of the free {boundaries},  we investigate next  their behavior as $t\to T$. Since $V^\mu(T,\cdot)=G$, there is no interface at $t=T$. Nevertheless, as it is shown next, there is a trace at $p=m$. Thus, defining $A_-^\mu(T)=A_+^\mu(T)=m$ gives {a} continuous extension
at $T$. %The definitions in \eqref{eq:boundary-definitions} apply only
%for $t<T$, since $\cE_T^\mu=\varnothing$.

\begin{proposition}
\label{prop:terminal-boundary-limits}
Assume \eqref{takis20}. Then,   as
$t\to  T$ and uniformly in $\mu\in\cM$, 
\[
\lim_{t\to T}\sup_{\mu\in\cM}
\max\{m-A_-^\mu(t),\,A_+^\mu(t)-m\}=0.
\]
\end{proposition}

\begin{proof}
Fix $\mu\in\cM$, $t<T$, $p\ne m$, write $P=P^{t,p;\mu}$, recall that 
the volatility of $P$ is bounded from {above} by $\overline c_b/4$, and, 
for any $\tau\in\mathcal T^\mu_{t,T}$, let
$H=\inf\{s\in[t,T]:P_s=m\}$, with $\inf\varnothing=\infty$. Then, 
It\^o's formula applied to $\tau$  yields 
\[
|p-m|^2\Pp(H\le\tau)\le\E|P_{\tau\wedge H}-p|^2
\le\frac{\overline c_b^2}{16}\E(\tau-t).
\]
Since $G$ is linear on each side of $m$, the optional sampling theorem gives
$\E G(P_{\tau\wedge H})=G(p)$, while the concavity of $G$ makes
$G(P_s)_{s\in [t,T]}$ a super-martingale. 
\vs
Hence, since $\{H\le\tau\}\in\mathcal F^\mu_\tau\cap\mathcal F^\mu_{H\wedge T}$
and $P_{H\wedge T}=m$ on this event, 
\[
\begin{aligned}
G(p)-\E G(P_\tau)
&=\E\bigl[\1_{\{H\le\tau\}}\bigl(G(m)-G(P_\tau)\bigr)\bigr]\\
&\le\E\bigl[\1_{\{H\le\tau\}}\bigl(G(m)-G(P_T)\bigr)\bigr]\\
&\le\frac{\ol a\overline c_b}{4}\sqrt{T-t}\,\Pp(H\le\tau),
\end{aligned}
\]
%Here $\{H\le\tau\}\in\mathcal F^\mu_\tau\cap\mathcal F^\mu_{H\wedge T}$
%and $P_{H\wedge T}=m$ on this event. 
with the last inequality following by
conditioning at $H\wedge T$, using the Lipschitz bound for $G$, the 
conditional It\^o isometry, and the Cauchy--Schwarz inequality.
\vs
Combining these estimates with $f^\mu\ge\underline c_f$ yields
\[
\E\!\left[\int_t^\tau f^\mu(s,P_s)\,\dd s+G(P_\tau)\right]-G(p)
\ge\left(\underline c_f-\frac{\ol a \overline c_b^3\sqrt{T-t}}{64|p-m|^2}\right)
\E(\tau-t).
\]
For every $0<\varepsilon<\min\{m,1-m\}$, the bracket is nonnegative, uniformly in $\mu$, 
whenever $|p-m|\ge\varepsilon$ and $t$ is sufficiently close to $T$.
\vs
Since $\tau=t$ is admissible, this proves
$V^\mu(t,p)=G(p)$ at all such points, and, thus Proposition~\ref{prop:two-boundary}
gives
$m-\varepsilon\le A_-^\mu(t)<m<A_+^\mu(t)\le m+\varepsilon$,
which proves the claim.

\end{proof}

\section{The existence of a mean field equilibrium}
\label{sec:mfg-existence}

Here, we present the proof of the main result of the paper. 
%We have now collected all the material needed to present the proof of the main result. 
For this we need to establish that, if $\mu_n\weakto \mu$ in $\cM$, then 
the best-response pair associated to  $\mu_n$, as defined in \eqref{eq:best-response-pair}, converges in probability to the best-response pair associated to $\mu$.

% cis continuous in probability if $\mu_n\weakto \mu$ in $\cM$. 
\vs

%Recall that, for  $\mu\in\mathcal C$, we denote by 
%$(\hat\tau_{0,p}^\mu,\hat d_{0,p}^\mu)$  the best-response pair defined
%in \eqref{eq:best-response-pair}.
%

The first step is to prove that the signal and belief processes have this property.

\begin{proposition}
\label{prop:canonical-response-continuity} 
Assume \eqref{takis20} and let $\mu^n, \mu \in {\mathcal M}$ be such that, as $n\to\oo$,   $\mu^n\weakto\mu$.  Consider the signal processes
\[
 \dd Y_n(t)=b\bigl(t,\Theta;F(t;(\mu^n)_{\le t})\bigr)\dd t+\dd W_t \quad \  Y_n(0)=0,\]
and
\[ \dd Y(t)=b\bigl(t,\Theta;F(t;\mu_{\le t})\bigr)\dd t+\dd W_t \ \ Y(0)=0,\]
and their respective belief processes $P_n:=P^{0,p_0;\mu^n}$ and
$P:=P^{0,p_0;\mu}$ 
with common initial belief $P_0\sim\eta$, and set
$\tau_n:=\hat\tau_{0,p_0}^{\mu^n}$ and
$\tau:=\hat\tau_{0,p_0}^{\mu}$.
Then, 
\begin{equation}
 \left(\tau_n,P_n(\tau_n)\right)
 \longrightarrow
 \left(\tau,P(\tau)\right)
 \quad\text{in probability},
 \label{eq:canonical-response-stability}
\nonumber
\end{equation}
and, consequently, $\Phi(\mu^n)\weakto\Phi(\mu)$ in {$\mathcal M$}.
\end{proposition}

\begin{proof}
The belief processes  $P_n$ and $P$ solve %(see {\eqref{vlakas0}} and {\eqref{0}}) 
%with dynamics 
\[
%\begin{split}
 \dd P_{n,t}
 ={\big(\Delta b(t;\mu^n_{\le t})\big)^2}P_{n,t}(1-P_{n,t})
   \bigl(\Theta-P_{n,t}\bigr)\dd t + \Delta b({t};\mu^n_{\le t})P_{n,t}(1-P_n(t))\dd W_t \ \ \ {P_{n,0}}=P_0, \]
   and
   \[  \dd P_t
 ={\big(\Delta b(t;\mu_{\le t})\big)^2}P_t(1-P_t)
   \bigl(\Theta-P_t\bigr)\dd t + \Delta b({t};\mu_{\le t})P_t(1-P_t)\dd W_t
  \ \ \ P_0=P_0.\]
%and similarly for $P$, with $\mu^n$ replaced by $\mu$.
%Both coefficients are uniformly Lipschitz in the state on $[0,1]$.

Combining  Lemma~\ref{lem:coefficient-stability},  the Burkholder-Davis-Grundy 
and Gr\"onwall's inequalities  and Proposition~\ref{prop:value-stability},  we find that, as $n\to \infty$, 
\[
 \E\Big[\sup_{t\le T}|P_n(t)-P(t)|^2\Big]\longrightarrow0 \ \ \text{and} \ \ \|\wh V^{\mu^n}-\wh V^\mu\|_\infty\to0.\]
Next, for  $t\in [0,T]$, we write %$\tau_n:=\hat\tau_{0,P_0}^{\mu^n}$, 
%$\tau:=\hat\tau_{0,P_0}^{\mu}$, 
{$X_n(t):=-\wh V^{\mu^n}(t,P_{n,t})$ and $X(t):=-\wh V^\mu(t,P_t)$}, and  denote by 
$\omega_{\wh V^\mu}$ the modulus of continuity of $\wh V^\mu$ on
$[0,T]\times[0,1]$,  which  exists since $\wh V^\mu$ is uniformly continuous
on this compact set.  
\vs
Then,
\begin{eqnarray}\label{eq:X_conv_prob}
   {\underset{s\in [0,T]} \sup |X_n(s)-X(s)|}
 \le {\|\wh V^{\mu^n}-\wh V^\mu\|_\infty}
    +\omega_{\wh V^\mu}\!\left(\|P_n-P\|_\infty\right)
 \underset{n\to \infty}\longrightarrow0
 \quad\text{in probability}.
\end{eqnarray}
Using the local uniform convergence on compact subsets of $\mathcal R\subset (0,T)$ of the $ A_\pm^{\mu^n}$'s to $A_\pm$, %$B^{\mu^n}$ to $A_-^\mu$ and $A_+^\mu$ respectively, 
%Let $\mathcal R\subset(0,T)$ be the full-measure set on which
%Theorem~\ref{thm:boundary-ae} gives
%\[
% A^{\mu^n}(s)\longrightarrow A_-^\mu(s),\qquad
% B^{\mu^n}(s)\longrightarrow A_+^\mu(s).
%\]
the preceding $L^2$ convergence of $P_n$ to $P$
and \eqref{eq:X_conv_prob}, we deduce that every
{subsequence} admits a further subsequence, which, for notational simplicity, we denote the same as the sequence,  and an event $\Omega_0$ with
$\Pp(\Omega_0)=1$ such that, {for every $\omega\in\Omega_0$,} %$[0,T)$, 
%  and as $n\to \infty$,

\begin{equation}\label{eq:coupled-pathwise-convergence}
 \underset{s\in [0,T]}\sup |P_{n,s}(\omega)-P_s(\omega)|
 + \underset{s\in [0,T]}\sup |X_n(s,\omega)-X(s,\omega)|\underset{n\to \infty} \longrightarrow0.
\end{equation}

Let $E_1$ be a probability-one event on which $P$ and
every $P_n$ have continuous paths; such an event exists because each
process is continuous almost surely and the sequence is countable.
Then, we replace $\Omega_0$ by $\Omega_0\cap E_1$, which still has probability
one, and we have that  the convergence in \eqref{eq:coupled-pathwise-convergence} and
path continuity hold on $\Omega_0\cap E_1$.

\vs

%Using the local uniform convergence on compact subsets of $\mathcal R\subset (0,T)$ of the $ A_\pm^{\mu^n}$'s to $A_\pm$, %$B^{\mu^n}$ to $A_-^\mu$ and $A_+^\mu$ respectively, 
%%Let $\mathcal R\subset(0,T)$ be the full-measure set on which
%%Theorem~\ref{thm:boundary-ae} gives
%%\[
%% A^{\mu^n}(s)\longrightarrow A_-^\mu(s),\qquad
%% B^{\mu^n}(s)\longrightarrow A_+^\mu(s).
%%\]
%%the synchronous path estimate and 
%and \eqref{eq:X_conv_prob}, we see that every
%{subsequence} admits a further subsequence, which, for notational simplicity, we denote the same as the sequence,  and an event $\Omega_0$ with
%$\Pp(\Omega_0)=1$ such that, {for every $\omega\in\Omega_0$,} %$[0,T)$, 
%%  and as $n\to \infty$,
%\[
% {\underset{s\in [0,T]}\sup |P_{n,s}(\omega)-P_s(\omega)|
% + \underset{s\in [0,T]}\sup |X_n(s,\omega)-X(s,\omega)|\underset{n\to \infty} \longrightarrow0.}
%\]
%{Before fixing a sample path, we need to intersect $\Omega_0$ with the
%probability-one events carrying the continuity of the posterior process, the conclusion of
%Proposition~\ref{lem:local-contact-Holder}, and the stopping-time law of the
%iterated logarithm \eqref{pll1}, the latter two on $\{0<\tau<T\}$.
%The resulting event, still denoted for simplicity by $\Omega_0$, has probability one.}
%\vs

Fix $\omega\in\Omega_0$ and $\delta>0$. If
$\tau(\omega)>\delta$, then 
\[c_\delta(\omega)
 :=\min_{0\le s\le\tau(\omega)-\delta}X(s,\omega)>0.\]
Indeed,  otherwise we would have  $X(s,\omega)=0$ for some
$s\le\tau(\omega)-\delta<\tau(\omega)$, which contradicts  the definition
of $\tau(\omega)$. 
\vs
Hence, for all sufficiently large $n$,
$X_n(s,\omega)\ge\frac{c_\delta(\omega)}2 $ for $ 0\le s\le\tau(\omega)-\delta$, and, thus,  $\tau_n(\omega)\ge\tau(\omega)-\delta,$
while, if 
$\tau(\omega)\le\delta$, we obviously have
$\tau_n(\omega)\ge0\ge\tau(\omega)-\delta$.
\vs
Letting
$\delta\to 0$ gives
\beq\label{eq:stopping-time-liminf}
\liminf_{n\to\infty}\tau_n(\omega)\ge\tau(\omega).
\eeq
Since $\omega\in\Omega_0$ was arbitrary and $\Pp(\Omega_0)=1$, the
inequality holds almost surely along the extracted subsequence.

%Fix $\omega\in\Omega_0$ and suppress it from the notation.  If $\tau>\delta$, continuity of $X$ and the definition of $\tau$ give $c_\delta:=\min_{0\le s\le\tau-\delta}X(s)>0$.  Hence, for all large $n$, $X_n\ge c_\delta/2$ on $[0,\tau-\delta]$, so $\tau_n\ge\tau-\delta$.  On $\{\tau\le\delta\}$, the inequality also holds because $\tau_n\ge0\ge\tau-\delta$.  Letting $\delta\downarrow0$ yields $\liminf_{n\to\infty}\tau_n\ge\tau$ almost surely along the extracted subsequence.
\vs
We now prove that $$\limsup_{n\to\infty}\tau_n\le\tau \ \ \text{a.s. on} \ \ 
\{0 \le \tau \le T\}.$$
\vs % in both cases.

We begin with the proof of the claim in $\{0< \tau <T\}.$

Fix $\delta>0$ and recall that, in view of  {Proposition~\ref{lem:local-contact-Holder}},   if $\tau \in (0,T)$ is the first contact point of 
the process $P$, 
%$(P_s)_{0\le s\le T}$  
and depending on whether the contact is on  $A_-^\mu$ or $A_+^\mu$, {the corresponding free boundary is $C^{0,1/2}$ on $(\tau-\overline r,\tau+\overline r)\cap(0,T)$ with constant $\overline C$, where $\overline C,\overline r>0$ may depend on $\mu$ and the realized contact point.}
% and the corresponding boundary has a
%local $C^{1/2}$ constant $C$.  
\vs

Using the uniform lower bound of the $a^\mu$'s and 
%In view of the strict positivity of the volatility away from $0$ and $1$,  
the law of  iterated logarithm for Brownian motion (see, for example, \cite[Theorem~2.9.23]{KaratzasShreve1991})  we deduce that 
\beq\label{pll1}
 \liminf_{h\downarrow0}\frac{P_{\tau+h}-P_\tau}{\sqrt h}=-\infty \ \ \ \text{and} \ \ \ 
 \limsup_{h\downarrow0}\frac{P_{\tau+h}-P_\tau}{\sqrt h}=\infty.
\eeq

Next, we show that \eqref{pll1} holds almost surely for a fixed $\mu$ and $\tau$. Indeed, we 
%To justify \eqref{pll1} as an almost-sure statement, for the fixed
%environment $\mu$ and stopping time $\tau$, 
apply the Dambis--Dubins--Schwarz
theorem and the law of the iterated logarithm to
$M_h:=P_{(\tau+h)\wedge T}-P_\tau$, which is a martingale with respect to
$(\mathcal F^\mu_{(\tau+h)\wedge T})_{h\ge0}$ with quadratic variation $\langle M\rangle$ satisfying, on
$\{0<\tau<T\}$, % its quadratic variation satisfies
\[\langle M\rangle_h/h\to(\Delta b^\mu(\tau))^2P_\tau^2(1-P_\tau)^2>0.\]

In view of  Proposition~\ref{lem:local-contact-Holder}, there is a
probability-one event $E_2$ on which, whenever $0<\tau<T$, first
contact occurs at one of the two boundaries and the contacted boundary
is locally $C^{1/2}$.
The above argument gives a probability-one event $E_3$ on which
both limits in \eqref{pll1} hold whenever $0<\tau<T$.
Replace $\Omega_0$ by $\Omega_0\cap E_2\cap E_3$, which still has
probability one. The convergence in \eqref{eq:coupled-pathwise-convergence}
and the $\liminf$ bound in \eqref{eq:stopping-time-liminf} remain valid
on this smaller event.
\vs
Fix $\omega\in\Omega_0\cap\{0<\tau<T\}$ and $\delta>0$. {For} notational simplicity in what follows we  suppress
$\omega$.
\vs

%In what follows, we need to examine separately the cases where the contact is either on the lower interface  $A^\mu_-$ or upper interface $A^\mu_+$. 
%\vs
%It follows that  \eqref{pll1}  holds on a probability-one event, which was included in
%$\Omega_0$ before fixing $\omega$.
%
%%\[\liminf_{h\downarrow0}\frac{W_{\tau+h}-W_\tau}{\sqrt h}=-\infty \ \ \text{and} \ \ 
%% \limsup_{h\downarrow0}\frac{W_{\tau+h}-W_\tau}{\sqrt h}=+\infty.\]
%%For the posterior process, applying DDS to the shifted martingale Pε+h↑Pε gives the same conclusions
%%on {ω < T, Pε → (0, 1)}, because
%\vs
In what follows, we examine separately the cases where the contact point  is either on the lower ({left}) free boundary   $A^\mu_-$ or upper ({right}) free boundary $A^\mu_+$. 
%Since the arguments are similar, we only provide the details for the former. 
\vs
If the former holds, using \eqref{pll1},  we can choose $h>0$ small enough (and  smaller than $\delta$ and $\overline r$ in Proposition~\ref{lem:local-contact-Holder}) to  obtain %at  the lower boundary $A_-^\mu$  and the upper boundary $A_+^\mu$ we have respctively 
\[P_{\tau+h}-P_\tau<-\overline C\sqrt h.\] %  \ \ \text{and} \ \ P(\tau+h)-P(\tau)>C\sqrt h.\]
It follows from the local $C^{1/2}-$continuity of the interfaces with constant $\overline  C$ that 
\[
 P_{\tau+h}<A_-^\mu(\tau)-\overline C\sqrt h\le A_-^\mu(\tau+h).\]
 % \ \ \text{and} \ \  P(\tau+h)>A_+^\mu(\tau)+C\sqrt h\ge A_+^\mu(\tau+h).
If the contact is on $A^\mu_+$, then, in view of \eqref{pll1}, we choose $h>0$ small enough (and  smaller than $\delta$ and $\overline r$ in {Proposition~\ref{lem:local-contact-Holder}}) so that 
\[P_{\tau+h}-P_\tau>\overline C\sqrt h.\]
Then, using the local $C^{1/2}-$continuity of the interfaces with constant $\overline  C$, we get 
\[
 P_{\tau+h}>A_+^\mu(\tau)+\overline C \sqrt h\ge A_+^\mu(\tau+h).
\]
%In either case, $h$ may be chosen smaller than $\delta$ and the radius $\overline r$
%in Lemma~\ref{lem:local-contact-H\"older}. 
%Then, the continuity of the free boundaries gives a nonempty open interval $J\subset (\tau, (\tau+\delta)\wedge T)$ on which 
%\vs
%
In summary, there exists  a nonempty compact set 
$J\subset\Big((\tau,(\tau+\delta)\wedge T)\cap \mathcal R\Big)$ on which
\[
\text{either}  \ \ P_s<A_-^\mu(s)\quad\text{or}\quad P_s>A_+^\mu(s),
\]
%Since $\mathcal R$ has full measure, choose $s\in\mathcal R\cap J$.
%S $P_n(s)\to P(s)$ by the definition of $\Omega_0$, and
while, as shown in the previous section, as $n\to \infty$ and uniformly on $J$, 
\[{A_-^{\mu^n}}\to A_-^\mu \ \ \text{and} \ \ A_+^{\mu^n}\to A_+^\mu.\]
%$B^{\mu^n}(s)\to A_+^\mu(s)$ by the definition of $\mathcal R$.
It then follows that, for all large $n$, % Therefore
\[
\text{either}  \ \  P_n(s)<{A_-^{\mu^n}(s)}\quad\text{or}\quad
 P_n(s)>{A_+^{\mu^n}(s)},
\]
and, thus, $\tau_n\le s<\tau+\delta$, and, in turn, 
$\limsup_{n\to\infty}\tau_n\le\tau$ on this event.
\vs

Since on $\{\tau=T\}$, it holds that $\tau_n\le T=\tau$, it only remains  to show  that on $\{\tau=0\}$ we also have %Hence
$\limsup_{n\to\infty}\tau_n\le\tau$.
%remains to be proved on
%$\{\tau=0\}$.  

\vs
Let $\mathcal R$ be the
full measure set from {Theorem~\ref{thm:boundary-ae}} and choose
$s_j\in\mathcal R$ with $s_j\downarrow0$. On $\{\tau=0\}$,
{Proposition~\ref{prop:two-boundary}} gives
\[
 p_0\in[0,{A_-^\mu(0)}]\cup[{A_+^\mu(0)},1].
\]
Since $\eta$ has a density, then either
$p_0<A_-^\mu(0)$ or $p_0>A_+^\mu(0)$ holds  almost surely. 
\vs
We only consider  the case $p_0<A_-^\mu(0)$, as the other one follows similarly.  
\vs
In view of the continuity of the belief process  and Proposition~\ref{lem:initial-piecewise}, 
there exists $r>0$ such that
$P_s<A_-^\mu(s)$ for $0\le s<r$, and, for every sufficiently large $j$,
$s_j<r$, and, furthermore, 
\[
 P_{ n,s_j}\to P_{s_j} \ \ \text{and} \ \ 
 A_-^{\mu_n}(s_j)\to A_-^\mu(s_j).
\]
Hence,  $ P_{ n,s_j}<{A_-^{\mu_n}(s_j)}$ for all large $n$, and, thus, 
$\tau_n\le s_j$.
\vs
It then follows that $\limsup_n\tau_n\le s_j$ for every sufficiently large $j$, and
letting $j\to\infty$  we deduce 
\[
 \limsup_{n\to\infty}\tau_n\le0=\tau.
\]
%The case $P_0>B^\mu(0)$ is identical, using $B^\mu$.
%
%If $\tau=0$, then $\wh V^\mu(0,p)=0$, which implies that $p\in [0, A_-^\mu(0)]\cup[A_+^\mu(0),1]$. Since $\eta$ is absolutely continuous
%with respect to the Lebesgue measure on $[0,1]$,  we only need to consider the case that $p\in [0, A_-^\mu(0))\cup(A_+^\mu(0),1]$. 
%\vs
%Since the arguments are similar, we assume next that, for some $\delta>0$,  $p>A_+^\mu(0) + \delta$. 
%
%\vs
%In view of the fact that, as $n\to \infty$,  $B^{\mu^n}(0)\to A_+^\mu(0)$, for sufficiently large $n$, we  $p>B^{\mu^n}(0)$.
%Thus, $\tau_n=0$, and the claim follows.
%\vs
%
%
%Then  \eqref{eq:canonical-response-stability}  follows from the fact that  $P_{n,{\tau_n}}\to P_\tau$ in probability follows from 
%\[
% |P_{n,{\tau_n}}-P_\tau|
% \le \|P_n-P\|_\infty+|P_{\tau_n}-P_\tau|.
%\]
%%which yields that  $P_{n,{\tau_n}}\to P_\tau$ in probability.
%% which,  together with $\tau_n\to\tau$, gives
%%\eqref{eq:canonical-response-stability}.
Next, we prove the convergence of the decision rules $D(P_{n,\tau_n})$.  
\vs
If $\tau<T$, 
Proposition~\ref{prop:two-boundary} gives
$P_\tau\in[0,A_-^\mu(\tau)]\cup[A_+^\mu(\tau),1]$, which does not contain
$m$.  If $\tau=T$, then $P_T$ has no atom at $m$, because the logit
diffusion has a transition density on $(0,1)$.  If $\tau=0$, then, as before, it suffices to consider the cases that $p$ is strictly  either below $A_-^\mu(0)$ or above $A_+^\mu(0)$ and, hence, away from $m$.
\vs
Thus,  the decision
rule is continuous at $P_\tau$ almost surely in all  cases, and, as $n\to \infty$,  $D(P_{n,\tau_n})\longrightarrow D(P_\tau)$ in probability. 
\vs  

Finally, for $h\in C([0,T])$ and $i\in\{0,1\}$, we have
\[
 \int_{[0,T]}h(t)\,{\nu_i^{\mu^n}}(\dd t)
 =\E\left[h(\hat\tau_{0,P_0}^{\mu^n})
   \1_{\{D(P_{n,\tau_n})=i\}}\right].
\]
The right-hand side converges to the analogous expression for $\mu$,
which proves $\Phi(\mu^n)\weakto\Phi(\mu)$.

\end{proof}

We conclude with the proof of the existence of the mean field equilibria.

\begin{theorem}%[Existence of a mean field equilibrium]
\label{thm:mfg-existence}
Assume \eqref{takis20}. Then, %~\ref{ass:coefficients} and~\ref{ass:obstacle},
there exists $\mu^*\in{\mathcal M}$ such that $\Phi(\mu^*)=\mu^*$.
\end{theorem}

\begin{proof}
The set {$\mathcal M$} is nonempty, convex and compact with respect to  the weak
topology, while  Proposition~\ref{prop:canonical-response-continuity}
shows that {$\Phi:\mathcal M\to\mathcal M$} is continuous. Then,  it follows from the 
Schauder--Tychonoff fixed-point theorem that there exists 
$\mu^*\in{\mathcal M}$ with $\Phi(\mu^*)=\mu^*$.
\end{proof}

\section{The  finite game}
\label{sec:N-player}

We consider the finite analogue of the information acquisition mean field game and examine the proximity of its solution to the ones in the continuum, constructed in the previous section. As mentioned in the introduction, we identify each decision maker (DM) with the pair of his prior belief and the related posterior belief process, generated by the signal whose dynamics depend on the stopping times, parameter assessment and costs of the rest of the populations. In this sense, there are no players as such in the standard way players are considered in mean field games. For this, we choose to use the terminology “decision maker” rather than a player, and we map any DM to the corresponding probabilities (prior and posterior). 
\vs
Fix $N\ge2$.  DM $j$ has initial belief $p^j\sim\eta$ and a private
{innovation} Brownian motion $Z^j$ {independent of $p^j$, with the pairs $(p^j,Z^j)$ independent across DMs}.  
\vs
Let
$\Omega_0=C_0([0,T];\R)$ with Wiener measure $\mathbb W$, and set 
\[
 \Pp^N:=(\eta\otimes\mathbb W)^{\otimes N}
 \quad\text{on }\quad ([0,1]\times\Omega_0)^N.
\]
For a fixed player $i$ and initial belief $p$, let $\Pp_{i,p}^N$ be the
same product law with $p^i=p$, and write $\E_{i,p}^N$ for expectation
under this law.
\vs
Let $\mu^*$ be a  mean field equilibrium  constructed in 
Theorem~\ref{thm:mfg-existence}, and, for each $j=1,\ldots,N$, consider 
%Every player other than a possible
%deviator uses the mean field 
best response strategy pair 
%$(\hat\tau_{0,p}^{\mu^*},\hat d_{0,p}^{\mu^*})$.  Thus, for every
%player $j$, define the strategy
\begin{equation}
 (\bar\tau^{N,j},\bar d^{N,j})
 :=(\hat\tau_{0,p^j}^{\mu^*},\hat d_{0,p^j}^{\mu^*}),
 \label{eq:mean field-induced-profile}
\end{equation}
with  the belief process driven by $Z^j$.  
\vs
Since
$\Phi(\mu^*)=\mu^*$, for every Borel $A\subset[0,T]$ and
$k=0,1$, we deduce that, for each $j=1,\ldots,N$,
\[
 \Pp^N(\bar\tau^{N,j}\in A,\ \bar d^{N,j}=k)=\mu^{*,k}(A).
\]
Therefore, all pairs $(\bar\tau^{N,j},\bar d^{N,j})$ are independent
and all  have the equilibrium  distribution $\mu^*$.
\vs

Next, we {fix} $i${. For} $j=1,\ldots, i-1, i+1, \ldots, N${,} we consider the pairs given by 
\eqref{eq:mean field-induced-profile} and, 
%Their realized stopping times and
%decisions determine the environment faced by player $i$ defined, 
for $k=0,1$, we introduce the measures 
\begin{equation}
 \widehat\mu^{-i,N,k}
 :=\frac1{N-1}\sum_{j\ne i}
   \1_{\{\bar d^{N,j}=k\}}\delta_{\bar\tau^{N,j}} 
 \ \ \text{and} \ \ 
 \widehat\mu^{-i,N}:=(\widehat\mu^{-i,N,0},\widehat\mu^{-i,N,1}).
 \label{eq:empirical-environment}
\nonumber
\end{equation}
Let $(\mathcal F_t^{N,i,p})$ be the right continuous completion under
$\Pp_{i,p}^N$ of the $\sigma-$ algebra
\[
 \sigma\!\left(p^i,Z_s^i,
   \widehat\mu^{-i,N,k}([0,s]):0\le s\le t,\ k=0,1\right).
\]
It then follows from the discussion in Section~$2$ that  the map 
%Since $H$ is continuous and $H(0)=0$,
$t\mapsto F(t;\widehat\mu^{-i,N}_{\le t})$ is continuous and adapted, and, 
hence,  predictable.  
%The other players' processes are
%independent of $(p^i,Z^i)$, so $Z^i$ remains a Brownian
%motion.
\vs
Next, we introduce the processes 
%For each $i=1,\ldots,N$ we define the empirical-environment posterior and its
%mean field comparison process, driven by the same $Z^i$, by
\begin{align}
 \dd P_t^{i,N,p}
 &=\Delta b(t;\widehat\mu^{-i,N}_{\le t})
   P_t^{i,N,p}(1-P_t^{i,N,p})\dd Z_t^i, \ \ \
 P_0^{i,N,p}=p,
 \label{eq:N-posterior}\\
 \dd P_t^{i,*,p}
 &=\Delta b(t;\mu^*_{\le t})
   P_t^{i,*,p}(1-P_t^{i,*,p})\dd Z_t^i,
\ \ \ P_0^{i,*,p}=p,
 \label{eq:mf-comparison-posterior}
\end{align}
which represent the posterior empirical environment  and its
mean field comparison process.
\vs
Since the coefficients above are predictable and uniformly Lipschitz in the state variable,
\eqref{eq:N-posterior} and \eqref{eq:mf-comparison-posterior} have unique strong $\mathcal F^{N,i,p}$-adapted solutions.
\vs
Let $\mathcal A_i^{N,p}$ be the pairs $(\tau,d)$, with  $\tau$ and $d$ being respectively  an
$[0,T]$-valued $\mathcal F^{N,i,p}$-stopping time and
$\mathcal F_\tau^{N,i,p}$-measurable.  
\vs

Then, in complete analogy to \eqref{000}, we consider the cost
\begin{equation}
 J_i^N(p;\tau,d)
 :=\E_{i,p}^N\left[
   \int_0^\tau f(s,P_s^{i,N,p};\widehat\mu^{-i,N}_{\le s})\dd s
   +\overline G(P_\tau^{i,N,p},d)\right].
 \label{eq:N-player-cost}
\nonumber
\end{equation}
Finally, under $\Pp_{i,p}^N$, we define 
\[
 \bar\tau^{N,i}
 =\inf\{t\in[0,T]:V^{\mu^*}(t,P_t^{i,*,p})
                         =G(P_t^{i,*,p})\}
 \ \ \text{and} \ \ 
 \bar d^{N,i}=D(P_{\bar\tau^{N,i}}^{i,*,p}),
\]
and easily deduce that  $(\bar\tau^{N,i},\bar d^{N,i})\in\mathcal A_i^{N,p}$.
\vs

The result about the finite size  approximation is stated and proved next.

\begin{theorem}%[$N$-player approximation]
\label{thm:epsilon-Nash}
Let $\mu^*$ be the mean field equilibrium constructed in 
Theorem~\ref{thm:mfg-existence}.  Then, there exists  a constant $C>0$, which is independent
of $N,i,$ and $p$,
such that, for every $p\in[0,1]$, $N\ge2$, $i=1,\ldots,N$,
and $(\tau,d)\in\mathcal A_i^{N,p}$,
\begin{equation}
 J_i^N(p;\bar\tau^{N,i},\bar d^{N,i})
 \le J_i^N(p;\tau,d)+\frac{C}{\sqrt{N-1}}.
 \label{eq:epsilon-Nash}
\nonumber
\end{equation}
Consequently, the mean field solution  is an
$\eps_N$-Nash equilibrium with $\eps_N=C/\sqrt{N-1}$.
\end{theorem}

\begin{proof}
Throughout the proof, $C$ denotes {a} generic constant independent of $N, p$ and $i=1,\ldots,N$.
%and of the deviation, whose value may change from line to line.  
\vs
For fixed $k=0,1$ and $j=1,\ldots,N$ such that $j\neq i$, we define
%first control the empirical mean field.  For fixed \(k=0,1\) and \(j\neq i\), define
\[
T_j^k :=
\begin{cases}
\bar\tau^{N,j}, & \text{if } \bar d^{N,j}=k,\\
T+1, & \text{if } \bar d^{N,j}\neq k.
\end{cases}
\]
It follows that  \((T_j^k)_{j\neq i}\) are i.i.d. and, for \(t\in[0,T]\),
\[
\widehat\mu_{-i,N,t}^{\,k}
=
\frac{1}{N-1}\sum_{j\neq i}\mathbf 1_{\{T_j^k\le t\}} \ \ \text{and} \ \
\mathbb P(T_j^k\le t)=\mu^{*,k}([0,t]).
\]
Thus, \(t\mapsto\widehat\mu_{-i,N,t}^{\,k}\) is the corresponding empirical distribution function, and  the
Dvoretzky--Kiefer--Wolfowitz inequality (see \cite{Massart1990}) yields
\begin{equation}
 \E_{i,p}^N\left[
  \sup_{t\le T}|\widehat\mu^{-i,N,k}([0,t])
                   -\mu^{*,k}([0,t])|^2\right]\le\frac1{N-1}.
 \label{eq:DKW-L2}
\end{equation}
Moreover  \eqref{eq:causal-smoothing} and Fubini's theorem imply that,
for every finite measure $\nu$,
\[
 F_k(t;\nu_{\le t})=\int_0^t\nu([0,s])\phi(t-s)\dd s.
\]
Therefore, 
\[
 \sup_{t\le T}|F_k(t;\widehat\mu^{-i,N,k}_{\le t})
                   -F_k(t;\mu^{*,k}_{\le t})|
 \le\sup_{t\le T}|\widehat\mu^{-i,N,k}([0,t])
                   -\mu^{*,k}([0,t])|,
\]
and, in turn,  \eqref{eq:DKW-L2} implies
\begin{equation}
 \E_{i,p}^N\left[
  \sup_{t\le T}|F(t;\widehat\mu^{-i,N}_{\le t})
                  -F(t;\mu^*_{\le t})|^2\right]\le\frac4{N-1}.
 \label{eq:aggregate-rate}
\end{equation}

For subprobability measures $\nu,\nu'$ on $[0,T]$, Kantorovich duality
and integration by parts give
\[
 \mathcal W_1^{\mathrm{ext}}(\nu,\nu')
 \le (2T+1)\sup_{t\in[0,T]}|\nu([0,t])-\nu'([0,t])|.
\]
Thus, \eqref{eq:DKW-L2} and the Lipschitz assumption \eqref{ass.b} imply
\[
\begin{aligned}
 &\E_{i,p}^N\left[\sup_{t\le T}
 \left|\Delta b(t;\widehat\mu^{-i,N}_{\le t})
       -\Delta b(t;\mu^*_{\le t})\right|^2\right]\\
 &\qquad\le C\E_{i,p}^N\left[
 d_{\mathcal M}(\widehat\mu^{-i,N},\mu^*)^2\right]
 \le\frac{C}{N-1}.
\end{aligned}
\]

\vs
It follows from the Lipschitz continuity of the coefficients of \eqref {eq:N-posterior} and \eqref{eq:mf-comparison-posterior}, 
%Next, we compare the two posterior processes $P^{i,N,p}$ and $P^{i,\ast,p}$.  Since
%the Lipschitz continuity $|x(1-x)-y(1-y)|\le|x-y|$ on $[0,1]$, 
the Burkholder--Davis--Grundy  inequality, and 
the  above estimate, that 
%Lipschitz continuity of $\Delta b$ with respect to its second argument , and
%\eqref{eq:aggregate-rate} give
\[
 \E_{i,p}^N\left[\sup_{r\le t}
 |P_r^{i,N,p}-P_r^{i,*,p}|^2\right]
 \le C\int_0^t\E_{i,p}^N\left[\sup_{u\le s}
 |P_u^{i,N,p}-P_u^{i,*,p}|^2\right]\dd s
   +\frac{C}{N-1}, 
\]
In turn, using Gr\"onwall's inequality, we find that, uniformly in $i$ and $p$,
\begin{equation}
 \E_{i,p}^N\left[\sup_{t\le T}
      |P_t^{i,N,p}-P_t^{i,*,p}|^2\right]\le\frac{C}{N-1}.
 \label{eq:posterior-N-rate}
\end{equation}

For any \((\tau,d)\in\mathcal A_i^{N,p}\), its cost under the
deterministic environment \(\mu^*\) is given by 
\[
 J_i^*(p;\tau,d)
 :=\E_{i,p}^N\left[
   \int_0^\tau f(s,P_s^{i,*,p};\mu^*_{\le s})\dd s
   +\overline G(P_\tau^{i,*,p},d)\right].
\]
The  Lipschitz bounds on $ f$ and $G$, together
with \eqref{eq:aggregate-rate}, \eqref{eq:posterior-N-rate}, and
Cauchy--Schwarz inequality, give  the uniform over $(\tau,d)\in\mathcal A_i^{N,p}$ estimate 
\begin{align}
 |J_i^N(p;\tau,d)-J_i^*(p;\tau,d)|
 &\le C\E_{i,p}^N\!\left[
   \sup_{t\le T}|P_t^{i,N,p}-P_t^{i,*,p}|
   +
   \sup_{t\le T}|F(t;\widehat\mu^{-i,N}_{\le t})
                    -F(t;\mu^*_{\le t})|\right]\nonumber \\
 &\le\frac{C}{\sqrt{N-1}}.
 \label{eq:uniform-cost-comparison}
\end{align}
%This estimate is uniform over $(\tau,d)\in\mathcal A_i^{N,p}$.

Finally, the information generated by the belief processes 
does not lower the value in the fixed environment $\mu^*$. Indeed,  we 
fix their 
sample paths,  and note that,  for almost every realization $\omega^{-i}$ if 
$\tau(\cdot,\omega^{-i})$ is a stopping time for the  private
filtration of $P^i$, then, for $k=0,1$ and $t\in[0,T]$, 
\[
 \{d(\cdot,\omega^{-i})=k\}
 \cap\{\tau(\cdot,\omega^{-i})\le t\}
 \in\sigma(p^i,Z_s^i:0\le s\le t).
\]
Therefore, for almost every \(\omega^{-i}\), fixing \(\omega^{-i}\) makes
\((\tau(\cdot,\omega^{-i}),d(\cdot,\omega^{-i}))\) admissible with
respect to this filtration. %player \(i\)'s private filtration. 
\vs
On the other hand,  the canonical pair is optimal in
the environment $\mu^*$, we deduce that 
\[
 J_i^*(p;\bar\tau^{N,i},\bar d^{N,i})=V^{\mu^*}(0,p)
 \le
 J_i^*(p;\tau(\cdot,\omega^{-i}),d(\cdot,\omega^{-i})).
\]
{Integrating over the  realizations of the other DMs with respect to
$(\eta\otimes\mathbb W)^{\otimes(N-1)}$ gives}
\[
 J_i^*(p;\bar\tau^{N,i},\bar d^{N,i})
 \le J_i^*(p;\tau,d).
\]
Combining the above  with \eqref{eq:uniform-cost-comparison} yields 
%first for the
%candidate pair and then for the deviation, gives, after enlarging, if necessary, $C$, 
\[
 J_i^N(p;\bar\tau^{N,i},\bar d^{N,i})
 \le J_i^N(p;\tau,d)+\frac{C}{\sqrt{N-1}}.
\]
and we conclude. %  \eqref{eq:epsilon-Nash} is proved.
%\vs
%For an arbitrary ex-ante
%type-dependent deviation, its regular conditional section given
%$p^i=p$ belongs to $\mathcal A_i^{N,p}$ for $\eta$-almost every $p$.
%Thus, applying the interim inequality to these sections and integrating against
%$\eta$ gives the ex-ante 
%$\eps_N$-Nash inequality.
%
\end{proof}

\bibliographystyle{plain}
\bibliography{bibsxz}

\end{document}